\documentclass[11pt,oneoside]{article}

 \usepackage{graphicx}
\usepackage{amsmath}
\usepackage{epstopdf} 
\usepackage{mathptmx}      
\usepackage{amsthm}
\usepackage{amsmath}
\usepackage{amsfonts}
\usepackage{amssymb}
\usepackage{mathrsfs}
\usepackage{amsfonts,color}
\usepackage{epstopdf}
\usepackage[citecolor=red]{hyperref}
\usepackage[citecolor=red]{hyperref}
\usepackage{cases}
\usepackage{booktabs}
\usepackage{multirow}
\numberwithin{equation}{section}
\numberwithin{table}{section}
\begin{document}

\title{A contribution to   condition numbers of the multidimensional total least squares problem with linear equality constraint}


\author{Qiaohua Liu\thanks{  Department of Mathematics, Shanghai University, Shanghai 200444, P.R. China.
              Supported by the National Natural Science Foundation of China under grant 11001167 (qhliu@shu.edu.cn).} \and Zhigang Jia\thanks{School of Mathematics and Statistics, Jiangsu Normal University,
Xuzhou 221116, P. R. China;
Research Institute of Mathematical Science, Jiangsu Normal University,
Xuzhou 221116, P. R. China.
Supported by the National Natural Science Foundation of China under grant 11771188 (zhgjia@jsnu.edu.cn).}\and
Yimin Wei\thanks{Corresponding author. School of Mathematical Sciences and Key Laboratory of Mathematics for Nonlinear Sciences,
Fudan University, Shanghai, 200433, P.R. China. Supported by the National Natural Science Foundation of China under grant 11771099 (ymwei@fudan.edu.cn).}}
\date{ }

\maketitle

\begin{abstract}
This paper is devoted to   condition numbers  of the  multidimensional total least squares  problem with linear equality constraint (TLSE).
Based on the perturbation theory of  invariant subspace, the TLSE problem is proved to be equivalent to a multidimensional unconstrained weighed  total least squares   problem in the limit sense. With a limit technique,
Kronecker-product-based formulae for normwise, mixed and componentwise condition numbers of  the minimum Frobenius norm TLSE solution are given.
 Compact upper bounds of these condition numbers are provided to  reduce the storage and computation cost.   All  expressions and upper bounds of these condition numbers  unify the ones  for the single-dimensional TLSE problem  and multidimensional total least squares problem.    Some numerical experiments are performed to illustrate our results.\\


\end{abstract}

\section{Introduction}

The multidimensional total least square (TLS) model, which arises in many data fitting and estimation problems,
finds a ``best'' fit to the  overdetermined system $Ax\approx B$,  where $A\in \mathbb{R}^{q\times n}(q>n)$   and $B\in \mathbb{R}^{q\times d}$   are contaminated by some noise.
It determines perturbations $E$ to the coefficient matrix $A$ and $F$ to the matrix $B$ measured by Frobenius norm  such that
 \begin{equation}
\min_{E,F} \|[E\quad F]\|_{F}, \quad \mbox{subject\ to}\quad (A+E)X=B+F. \label{1.1}
\end{equation}
After the minimizer  $[\hat E\quad \hat F]$ is found such that
the corrected system $(A+\hat E)X=B+\hat F$ is consistent, the corresponding solution $X$ is called the TLS solution.
The TLS   model, was originally proposed in 1901 for data fitting problem \cite{pe}, but  hasn't caught much attention for a long time.
In 1980, Golub and Van Loan \cite{gv} introduced this model into the numerical linear algebra area.  Since then,  it has been attracting more and more attention and now the TLS model is widely applied in a broad class of scientific disciplines such
as system  identification \cite{ld}, image  processing \cite{nbk,npp},  speech and audio processing \cite{hv,lmv}, etc.  An overview of applications, theory, and computational methods
of the TLS problem, we refer to   \cite{gv,mv, val,vh1, xxw2}.

An extension of TLS model is the following multidimensional TLS problem with equality constraint (TLSE):
\begin{equation}
\min_{E,F} \|[E\quad F]\|_{F}, \quad \mbox{subject\ to}\quad (A+E)X=B+F,\quad CX=D, \label{1.2}
\end{equation}
where $D\in {\mathbb R}^{n\times d}$ and $C\in {\mathbb R}^{p\times n}$ is of full row rank.  When $d=1$, it reduces to the single-dimensional TLSE problem, which was first presented by   Dowling,  Degroat, and Linebarger \cite{ddl} in 1992, where a stable algorithm  on the  basis of QR and singular value decomposition (SVD) matrix factorizations were proposed. Further investigations on the single-dimensional TLSE  were performed in  \cite{sc}, where   iteration methods  were derived based on the Euler-Lagrange theorem.  Recently, Liu et al. \cite{liu} investigated   uniqueness conditions of the single-dimensional TLSE solution and interpreted the  solution as an approximation of the solution to an unconstrained weighted TLS problem (WTLS), with a large weight
 assigned on the constraint, based on which a QR-based inverse iteration method was presented.

The sensitive analysis and the condition number of a problem are  vital in numerical analysis, since the  condition number measures
the worst-case sensitivity of its solution to small perturbations in the input data. Combined with backward error
estimate,   an approximate upper bound can be derived for the forward error.

When $C, D$ are zero matrices and $d=1$, the TLSE problem becomes the standard single dimensional TLS problem, whose first order perturbation analysis and condition numbers  have been widely studied \cite{bg,ds,dwx,jl,lj,xxw,zlw}.   The condition number of the truncated TLS solution of an ill-conditioned TLS problem was studied by  Gratton, Titley-Peloquin,
and Ilunga \cite{gti}, Meng, Diao and Bai \cite{mdb}.  By  making use of  the perturbation results in \cite{bg,jl,lj}, and  the close relation of
 the single dimensional TLSE to an unconstrained  weighted TLS problems, Liu and Jia \cite{liu2} derived closed formulae  for condition numbers of the single dimensional TLSE problem.
  Further perturbation results were given in \cite{liu1}, which provides perturbation analysis and  tighter bounds for the forward error of the solution, when the perturbation in   input data are of different magnitude. The condition numbers and perturbation results in \cite{liu1,liu2} unify those for standard TLS problem \cite{jl,lj,xxw,zlw}.   When $C, D$ are nonzero matrices and $d=1$, under some condition (see (\ref{2.002}) with $\widetilde\sigma_{n-p+1}=0$), the TLSE solution reduces to a solution to the least squares problem with equality constraint (LSE), whose perturbation results were   studied in  \cite{ch,diao3,wei2}, that are also unified by the ones \cite{liu1} for the TLSE  problem.

  When $C, D$ are zero matrices and $d>1$, the TLSE problem becomes the multidimensional TLS problem.  In \cite{zmw}, Zheng, Meng and Wei studied the explicit formulae for the condition numbers of the minimum Frobenius norm    TLS solution.  The condition numbers of  multidimensional TLS problem with more than one solution were further studied in \cite{mzw} by Meng, Zheng and Wei.

To the best of our knowledge,    condition numbers of the multidimensional TLSE problem haven't been addressed in literature. In this paper, we aim to study this issue.  With the invariant subspace perturbation theorem, we prove that it is equivalent to a multidimensional  weighted TLS problem, with a large weight assigned on the constraint. By making use of the perturbation estimates in   \cite{zmw} for the multidimensional TLS problem,   we establish the first order perturbation estimates  of the minimum Frobenius norm   TLSE solution based on a limit technique, from which Kronecker-product-based  normwise, mixed and componentwise condition numbers formulae are derived. In order to reduce the storage and computation cost in these Kronecker-product-based formulae, compact upper bounds
of these condition numbers are given. The newly derived results unify those for multidimensional TLS and single-dimensional TLSE problems. Numerical examples are provided to show their tightness.

 Throughout the paper,  $\|\cdot\|_2$ denotes the Euclidean vector or matrix  norm,  $I_n$, $0_n$, $0_{m\times n}$ denote  the $n\times n$ identity matrix, $n\times n$ zero matrix, and  $m\times n$ zero matrix, respectively. If subscripts are ignored, the sizes of identity and zero matrices are suitable with context.
   For a matrix $M\in {\mathbb R}^{m\times n}$, $M^T$, $M^\dag$, ${\cal R}(M)$, $\sigma_i(M)$($\sigma_{\rm min}(M)$), $\|M\|_{\rm max}$  denote the transpose, Moore-Penrose inverse, the column range space, the $i$-th largest (the smallest) singular value, the maximal absolute value of elements of $M$, respectively.
${\rm vec}(M)$ is an operator, which stacks the columns of $M$ one  underneath the other.  For
  matrices $A\in {\mathbb R}^{m\times n}$ and $B\in {\mathbb R}^{s\times t}$,
 the Kronecker product   of $A, B$ is defined by $A\otimes B=[a_{ij}B]$ and  its property is listed as follows \cite{gr,ls}:
$$
\begin{array}l
{\rm vec}(AXB)=(B^{T}\otimes A){\rm vec}(X),\quad (A\otimes B)(C\otimes D)=(AC)\otimes (BD),\\
(A\otimes B)^T=A^T\otimes B^T,\quad (A\otimes B)^\dag=A^\dag\otimes B^\dag,\quad \|A\otimes B\|_2=\|A\|_2\|B\|_2,\\
 {\rm vec}(A^T)=\Pi_{(m,n)}{\rm vec}(A),\quad  \Pi_{(s,m)}(A\otimes B)=(B\otimes A)\Pi_{(t,n)},
\end{array}
$$
where $X\in {\mathbb R}^{n\times s}$, $C\in {\mathbb R}^{n\times k}, D\in {\mathbb R}^{t\times r}$ and $\Pi_{(m,n)}$ is an $mn\times mn$  vec-permutation matrix taking the form $\Pi_{(m,n)}=\sum_{i=1}^m\sum_{j=1}^nE_{ij}\otimes E_{ij}^T$, in which $E_{ij}\in {\mathbb R}^{m\times n}$ has an entry in  position $(i,j)$ and all other entries are zero.

\section{Preliminaries}

\indent In this section we first recall some well known results about multidimensional TLS and single-dimensional TLSE problems, after which we give solvability conditions and explicit form for the multidimensional TLSE solution.

\subsection{The first order perturbation estimate for   multidimensional TLS problems}

Let $L\in {\mathbb R}^{m\times n}$, $H\in {\mathbb R}^{m\times   d}$ ($m\ge n+d$), the multidimensional TLS problem is defined by
\begin{equation}
\min_{E, F}\|[E\quad F]\|_F,\qquad {s.t.}\qquad (L+E)X=H+F.\label{2.1}
\end{equation}
Following  \cite{gv}, the TLS problem (\ref{2.1}) may have no solutions. In order to broad its scope of applications,
the generic and nongeneric conditions for  TLS solutions were further  studied by  Van Huffel and Vandewalle \cite{vv,vv2}.
In 1992, Wei \cite{wei} redefined the conditions (see Eq. (\ref{c1})) to make the TLS problem (\ref{2.1})  meaningful in any situation.
The condition in (\ref{c1}) includes those in \cite{gv,vv,vv2} as   special cases.

SVD is a useful tool to characterize  the TLS solution. If the skinny SVD \cite[Chapter 2.4]{gv2} of $[L\quad H]$ is given by
\begin{equation}
[L\quad H]=U\Sigma V^T,\quad \Sigma=\mbox{diag}(\sigma_1,\sigma_2,\cdots,\sigma_{n+d})\in {\mathbb R}^{(n+d)\times (n+d)},\label{2.2}
\end{equation}
where $\sigma_i=\sigma_{i}([L\quad H])$ and $\sigma_1\geq \sigma_2\geq \cdots\geq \sigma_{n+d}\ge 0,$  $U\in {\mathbb R}^{m\times (n+d)}$ and $V\in {\mathbb R}^{(n+d)\times (n+d)}$ have orthonormal columns.  For an integer $t\in [0,n]$, partition
$$
\begin{array}l
V=\begin{array}l n\\ d\end{array}\left[\begin{array}{cc} V_{11}(t) & \quad V_{12}(t)\\ V_{21}(t) &\quad V_{22}(t)\end{array}\right].\\
\qquad\quad\quad~~~~ t \qquad~~ n+d-t
\end{array}
$$
For simplicity, we denote $V_{ij}=V_{ij}(k)$ for $i,j=1,2$. If
\begin{equation}\sigma_t>\sigma_{t+1},\quad  {\rm rank}(V_{22})=d,\label{c1}\end{equation}
holds, a solution to the consistent linear system
$\widehat LX=\widehat H$
is defined as a TLS solution to the linear approximation equation $LX\approx H$, where $\widehat L=U_1\Sigma_1V_{11}^T$ and $\widehat H=U_1\Sigma_1V_{21}^T$ with $U_1, V_1$ being, respectively, the first $t$ columns of $U=[U_1\quad U_2]$ and $V=[V_1\quad  V_2]$. The diagonal matrices  $\Sigma_1={\rm diag}(\sigma_1,\sigma_2,\cdots,\sigma_{t})$ and $\Sigma_2={\rm diag}(\sigma_{t+1},\sigma_2,\cdots,\sigma_{n+d})$. Among all TLS solutions, the minimum Frobenius norm solution to the compatible system  is given by $X_t=-V_{12}V_{22}^\dag$.

In \cite{zmw}, Zheng, Meng and Wei defined the mapping $\phi$: ${\mathbb R}^{m(n+d)}\rightarrow {\mathbb R}^{nd}$ by $\phi(c)={\rm vec}(X_t)$ for $c={\rm vec}([L\quad H])$ and
provided the first order perturbation analysis of  $\phi(c)$ as
\begin{equation}
\begin{array}{rl}
{\rm vec}(\Delta X_t)&=\phi'(c){\rm vec}([\Delta L\quad \Delta H])+{\cal O}(\|\Delta L\|_F^2+ \|\Delta H\|_F^2)\\
&=(H_1+H_2)DZ{\rm vec}([\Delta L\quad \Delta H])+{\cal O}(\|\Delta L\|_F^2+ \|\Delta H\|_F^2),
\end{array}\label{2.3}
\end{equation}
where
$$
\begin{array}l
H_1=\Big((  V_{22} V_{22}^T)^{-1} V_{21}\otimes   (V_{12}F_{ V_{22}}) \Big),\qquad
H_2=\Big( V_{22}^{\dag T}\otimes   V_{11}^{\dag T}\Big)\Pi_{(n+d-t,t)},\\
D=( \Sigma_1^2\otimes I_{n+d-t}-I_t\otimes ( \Sigma_2^T  \Sigma_2))^{-1}\big[ I_t\otimes  \Sigma_2^T\qquad \Sigma_1\otimes I_{n+d-t}\big],\\
Z=\left[\begin{array}c V_1^T\otimes U_2^T\\ \Pi_{(t,n+d-t)}(V_2^T\otimes U_1^T)\end{array}\right],
\end{array}
$$
in which
 $\Pi_{(n+d-t,t)}$ is a  vec-permutation matrix, $F_{  V_{22}}=I-  V_{22}^\dag   V_{22}$. From this result,  the absolute normwise condition number
  $\kappa_{\rm tls}^{\rm abs}(X_t,L,H)$ satisfies
  \begin{equation}
 \kappa_{\rm tls}^{\rm abs}(X_t,L,H)= \|(H_1+H_2)D\|_2\le (1+\|X_t\|_2^2){\displaystyle \sqrt{\sigma_t^2+\sigma_{t+1}^2}\over\displaystyle \sigma_t^2-\sigma_{t+1}^2},\label{bd1}
  \end{equation}
where the upper bound  is proved to be optimal and is attainable for some specific matrices.  In particular, for $t=n$,
  \begin{equation}
  {\displaystyle \sqrt{ \sigma_n^2+ \sigma_{n+1}^2}\over\displaystyle \|V_{11}\|_2\|V_{22}\|_2(\sigma_n^2- \sigma_{n+1}^2)}\le
  \kappa_{\rm tls}^{\rm abs}(X_n,L,H)\le (1+\|X_t\|_2^2){\displaystyle \sqrt{\sigma_n^2+\sigma_{n+1}^2}\over\displaystyle \sigma_n^2-\sigma_{n+1}^2}.\label{bd2}
  \end{equation}

\subsection{Solvability conditions and explicit solution of multi-dimensional TLSE problem}

For the multidimensional TLSE problem (\ref{1.2}), denote $\widetilde A=[A~~ B], \widetilde C=[C~~D]$, and assume that the QR factorization of  $\widetilde C^T$ takes the form:
 \begin{equation}
 \widetilde C^T=\widetilde Q\Big[{\widetilde R_1\atop 0}\Big],\qquad \widetilde Q=[\widetilde Q_1\quad \widetilde Q_2],\label{2.03}
 \end{equation}
 in which  $\widetilde Q_1\in {\mathbb R}^{(n+d)\times p}$, $\widetilde Q_2\in {\mathbb R}^{(n+d)\times (n+d-p)}$. Let  the skinny SVD of $\widetilde A\widetilde Q_2$ as
 \begin{equation}
 \widetilde A\widetilde Q_2=\widetilde U\widetilde\Sigma\widetilde V^T=[\widetilde U_1\quad \widetilde U_2]\left[\begin{array}{cc}
 \widetilde \Sigma_1&0\\0&\widetilde \Sigma_2\end{array}\right][\widetilde V_1\quad \widetilde V_2]^T,\label{2.4}
 \end{equation}
 where the matrices $\widetilde U_1, \widetilde U_2$ have, respectively, $k, (n+d-k-p)$ orthonormal columns.
  $\widetilde V_1$   is the submatrices of the $(n+d-p)\times (n+d-p)$ orthogonal matrix  $\widetilde V$ by  taking their first $k$ columns. The diagonal matrices  $\widetilde \Sigma={\rm diag}(\widetilde\sigma_1, \widetilde\sigma_2,$ $\ldots, \widetilde\sigma_{n+d-p})$,
 $\widetilde \Sigma_1={\rm diag}(\widetilde\sigma_1,$ $\widetilde\sigma_2,\ldots, \widetilde\sigma_{k})$,  $\widetilde \Sigma_2={\rm diag}(\widetilde\sigma_{k+1},\widetilde\sigma_{k+2},\ldots, \widetilde\sigma_{n+d-p})$, in which  $0\le k\le n-p$ is an integer such that
 \begin{equation}
 {\cal C}(k):\qquad \widetilde\sigma_1\ge  \widetilde\sigma_2\ge \ldots\ge  \widetilde\sigma_{k}>\widetilde\sigma_{k+1}\ge \ldots\ge\widetilde\sigma_{n+d-p}.\label{2.5}
 \end{equation}

For $d=1$, $x_C=C^\dag D$, and $ r_C=Ax_C-B$, in \cite{liu}, Liu et al. proved that  if
the orthonormal basis of null space of $\widetilde C$  is chosen as
\begin{equation}
\widetilde Q_2=\left[\begin{array}{cc} Q_2&\beta^{-1}x_C\\0&-\beta^{-1}\end{array}\right],\quad \beta=\big(1+\|x_C\|_2^2\big)^{1/2},\label{2.001}
\end{equation}
in which $Q_2$ is the orthonormal basis of the null space of $C$, then under the condition
\begin{equation}
 \sigma_{n-p}(AQ_2)>\sigma_{n-p+1}([AQ_2\quad \beta^{-1}r_C])=\sigma_{n-p+1}(\widetilde A\widetilde Q_2)=\widetilde\sigma_{n-p+1},\label{2a}
\end{equation}
the  TLSE solution $x_n$ is unique and can be expressed by
 \begin{equation}
 x_n=x_C-{\cal K}A^Tr_C,\quad for\quad {\cal K}=Q_2(Q_2^TA^TAQ_2-\widetilde\sigma_{n-p+1}^2I_{n-p})^{-1}Q_2^T.\label{2.002}
 \end{equation}

In the following theorem, we give the solvability conditions and  explicit form of the solution to the multidimensional case.

{\it {\bf Theorem 2.1} With the notation in (\ref{2.03})--(\ref{2.5}), let $t=p+k$ and
 $\overline V=\widetilde Q_2\widetilde V_2$ have the partition
\begin{equation}
\begin{array}{rl}
\overline V=\widetilde Q_2\widetilde V=&[\overline V_1\quad \overline V_2]=\begin{array}l n\\ d\end{array}\left[\begin{array}{cc}
\overline V_{11}&\quad \overline V_{12}\\ \overline V_{21}&\quad\overline V_{22}\end{array}\right].\\
&~ k~~~~n+d-t~~~ ~~~~~~~~k~\quad~n+d-t\\
\end{array}\label{2.6}
\end{equation}
  If for $k=n-p$,  the condition ${\cal C}(k)$ holds such that $\overline V_{22}$ is nonsingular, then the unique  TLSE solution is determined by $X_{n}=-\overline V_{12}\overline V_{22}^{-1}$, which is also the solution the consistent linear system
\begin{equation}
 \widehat AX=\widehat B,  \quad\mbox{subject to}\quad   CX=D,\label{2.7}
\end{equation}
  where
 \begin{equation}\widehat A=\widetilde U_1\widetilde\Sigma_1\overline V_{11}^T,\qquad
\widehat B=\widetilde U_1\widetilde\Sigma_1\overline V_{21}^T. \label{2.14}
\end{equation}
  }

{\it Proof.} Let $\widetilde X=\left[X^T~ ~-I_d\right]^T$.
Notice that  the constraint $CX=D$ requires   $\widetilde C\widetilde X=0$, therefore
 $\widetilde X$ lies in the null space of $\widetilde C$ spanned by $\widetilde Q_2$. Denote $\widetilde X=\widetilde Q_2Z$, and write
$\widetilde A=\widetilde A\widetilde Q_1\widetilde Q_1^T+\widetilde A\widetilde Q_2\widetilde Q_2^T,$  $\widetilde E=[E\quad F]$, (\ref{1.2}) becomes
\begin{equation}
\min \|[\widetilde E\widetilde Q_1\quad  \widetilde E\widetilde Q_2]\|_F, \quad s.t.\quad (\widetilde A\widetilde Q_2+\widetilde E\widetilde Q_2)Z=0,\label{2.8}
\end{equation}
where the restriction only imposed on $E\widetilde Q_2$ means that we can choose optimal $\widetilde E_*$ such that $\widetilde E_*\widetilde Q_1=0$ and $\widetilde A\widetilde Q_2+\widetilde E_*\widetilde Q_2$ has a null space with    dimension  no less than $d$.

Note that the condition $\widetilde E_*\widetilde Q_1=0$ means there exists a matrix $Y$ such that $\widetilde E_*^T=\widetilde Q_2Y^T$, and (\ref{2.8}) becomes
$$
\min\limits_{{\rm rank}(\widetilde A\widetilde Q_2+Y)\le n-p} \|Y\|_F, \quad s.t.\quad (\widetilde A\widetilde Q_2+Y)Z=0.
$$
 According to (\ref{2.4}) and the well-known Eckart-Young  theorem \cite[Theorem 2.4.8]{gv2} for the best rank-$(n-p)$ matrix approximation,  the optimal $Y_*$ satisfies
 $Y_*=-\widetilde U_2\widetilde\Sigma_2\widetilde V_2^T,$
and  for this optimal  error matrix $\widetilde E_*=Y_*\widetilde Q_2$,  the corrected system becomes
$$
(\widetilde A\widetilde Q_2-\widetilde U_2\widetilde \Sigma_2\widetilde V_2^T)Z=0,\quad or\quad \widetilde U_1\widetilde\Sigma_1(\widetilde Q_2\widetilde V_1)^T\widetilde X=0.
$$
Recall that $\widetilde X\in {\cal R}(\widetilde Q_2)$, therefore
   $\widetilde X$ lies in the range of $\overline V_2=\widetilde Q_2\widetilde V_2$,
    i.e., there exists an  $d\times d$  matrix $G$ such that
\begin{equation}
\left[{X\atop -I_d}\right]=\Big[{\overline V_{12}\atop \overline V_{22}}\Big]G,\label{2.15}
\end{equation}
from which we obtain  $G=-\overline V_{22}^{-1}$ and the unique solution is given by $X_n=-\overline V_{12}\overline V_{22}^{-1}$. \qed\\

{\bf Remark 2.1} If the condition ${\cal C}(k)$ only  holds for $k<n-p$  such that the right bottom partition $\overline V_{22}$ in (\ref{2.6}) is of full row-rank, we define a solution to
the linear system  (\ref{2.7})-(\ref{2.14})   as  a  TLSE solution.
In this case, $CX=D$ requires that $\Big[{X\atop -I_d}\Big]\in {\cal R}(\widetilde Q_2)={\cal R}(\overline V)$, and at the same time we notice that
$$[~ \widehat A\quad \widehat B~ ]\Big[{X\atop -I_d}\Big]=\widetilde U\widetilde\Sigma \overline V_1^T\Big[{X\atop -I_d}\Big]=0,$$
therefore $\Big[{X\atop -I_d}\Big]$ lies in the range space ${\cal R}(\overline V_2)$ and there exists an matrix $G$ such that (\ref{2.15}) holds.
From the relation $\overline V_{22}G=-I_d$, we conclude that
$G=-\overline V_{22}^\dag+PK$ for $P=I_{n+d-t}-\overline V_{22}^\dag\overline V_{22}$ and an arbitrary $(n+d-t)\times d$ matrix $K$.
Therefore  any TLSE solution $X$ has the form
$$X=-\overline V_{12}\overline V_{22}^\dag +\overline V_{12}PK,$$ in which
$$
\Big(\overline V_{12}\overline V_{22}^\dag\Big)^T\overline V_{12}P=\overline V_{22}^{\dag^T}\overline V_{12}^T\overline V_{12}P=\overline V_{22}^{\dag^T}(I-\overline V_{22}^T\overline V_{22})P=0,
$$
and    $X_{t}=-\overline V_{12}\overline V_{22}^\dag$ is the minimum Frobenius norm solution among all TLSE solutions.\\

\section{Close relation of TLSE to an unconstrained weighted TLS problem}

When $d=1$, Liu et al. \cite{liu}  proved that under the genericity condition (\ref{2a}),   the  unique solution of the single dimensional TLSE problem  can be interpreted as an approximation of the solution to an unconstrained weighted TLS   problem,  by assigning a large weight on the constraint.

Similar conclusions can be drawn for the multidimensional case. However,
in proving this assertion, we cann't mimic the technique in \cite{liu}, since some singular values of $\widetilde A\widetilde Q_2$ characterized by $\widetilde \Sigma_2$  might be multiple in the multidimensional case, and the associated singular vectors are not uniquely determined.  We need to  generalize Stewart's result \cite{st}  about
   the asymptotic behavior for the  scaled SVD of $X_\epsilon=[X_1\quad \epsilon X_2]$, based on the  following perturbation theorem for invariant subspaces.

{\bf Lemma 3.1} \cite[Chp. V, Thm 2.7]{sun} Let $[Z_1\quad Y_2]\in {\mathbb R}^{n\times n}$ be an orthogonal matrix and ${\cal R}(Z_1)$ is a $k$-dimensional simple invariant subspace of $n\times n$ matrix $C$ such that
$$
[Z_1\quad Y_2]^TC[Z_1\quad Y_2]=\left[\begin{array}{cc} L_1&H\\ 0&L_2\end{array}\right],
$$
where $L_1$ and $L_2$ have no common eigenvalues, $Y_2^TCZ_1=0$ ( Here ${\cal R}(Z_1), {\cal R}(Y_2)$ are called the right and left invariant subspace of $C$). Given a perturbation $E$, let
$$
[Z_1\quad Y_2]^TE[Z_1\quad Y_2]=\left[\begin{array}{cc} E_{11}&E_{12}\\
E_{21}&E_{22}\end{array}\right].
$$
Then for perturbations  $\|E\|_2$ small enough, there is a unique matrix $P$ such that the columns of
\begin{equation}
\tilde Z_1=(Z_1+Y_2P)(I+P^TP)^{-{1\over 2}},\qquad
\tilde Y_2=(Y_2-Z_1P^T)(I+PP^T)^{-{1\over 2}} \label{3.01}
\end{equation}
form orthonormal bases for simple right and left invariant subspaces of $\tilde C=C+E$. The representation of $\tilde C$, i.e.,
  $\tilde C\tilde Z_1=\tilde Z_1\tilde L_1$, $\tilde C\tilde Y_2=\tilde Y_2\tilde L_2$ with respect to
$\tilde Z_1$, $\tilde Y_2$  are given by
\begin{equation}
\begin{array}l
\tilde L_1=(I+P^TP)^{1\over 2}[L_1+E_{11}+(H+E_{12})P](I+P^TP)^{-{1\over 2}},\\
\tilde L_2=(I+PP^T)^{-{1\over 2}}[L_2+E_{22}-P(H+E_{12})](I+PP^T)^{1\over 2}.
\end{array}\label{3.02}
\end{equation}

{\bf Lemma 3.2} Let $\epsilon>0$ be a small parameter,  $X=[X_1\quad X_2]\in {\mathbb R}^{m\times n}$ with $X_1\in {\mathbb R}^{m\times k}$  being of full column-rank. Denote $X_\epsilon=[X_1\quad \epsilon X_2]$, $\overline X_2=X_2-X_1B$  with $B=X_1^\dag X_2$.
Assume that $X_1=U_1S_1V_1^T$, ${\overline X}_2={\overline U}_2{\overline S}_2{\overline V}_2^T$ and $X_\epsilon=U_\epsilon S_\epsilon V_\epsilon^T$ are the skinny SVDs of $X_1$, $\overline X_2$ and $X_\epsilon$, respectively, then
$$
\begin{array}l
S_\epsilon={\rm diag}(S_1+{\cal O}(\epsilon^2), \epsilon \overline S_2+{\cal O}(\epsilon^3)),\\
  U_\epsilon=\left[ U_1+{\cal O}(\epsilon^2)  \quad \overline U_2+{\cal O}(\epsilon^2)\right],\quad
  V_\epsilon=\left[\begin{array}{cc}
V_1+{\cal O}(\epsilon^2)&-\epsilon B\overline V_2+{\cal O}(\epsilon^3)\\
\epsilon B^TV_1+{\cal O}(\epsilon^3)&\overline V_2+{\cal O}(\epsilon^2)\end{array}\right].
  \end{array}
$$

{\it Proof}. Let $G=[X_1\quad 0_{m\times (n-k)}]^T[X_1\quad  0_{m\times (n-k)}]$ and
$$
G_\epsilon=X_\epsilon^TX_\epsilon=G+\left[\begin{array}{cc} 0& \epsilon X_1^TX_2\\ \epsilon X_2^TX_1& \epsilon^2X_2^TX_2\end{array}\right]=:G+E
$$
be the perturbed version of $G$.
Notice that
$$
\begin{array}{rl}
[Z_1\quad  Y_2]&=\begin{array}c k\\ n-k\end{array}\left[\begin{array}{cc} V_1&0 \\0& \overline V_2\end{array}\right]\\
&\qquad\qquad~~~~ k\quad n-k
\end{array}
$$
has  orthonormal columns and   $Z_1$, $Y_2$ form orthonormal bases of simple invariant subspace of $G$  such that the representations of $G$ with respect to $Z_1, Y_2$ are
\begin{equation}
GZ_1=Z_1L_1, \qquad GY_2=Y_2L_2,\qquad \mbox{for}\quad L_1=S_1^TS_1,\quad  L_2=0_{n-k}.\label{c2}
\end{equation}

By Lemma 3.1, there exists an $(n-k)\times k$ matrix $P$ such that  $\tilde Z_1, \tilde Y_2$ with   structure  (\ref{3.01})  form the orthonormal bases of right and left invariant subspaces of
$G_\epsilon$, respectively.  Substituting
the matrices $Z_1, Y_2$ and formula (\ref{3.01}) into the relation   $\tilde Y_2^TG_\epsilon \tilde Z_1=0$, one can derive that $(Y_2-Z_1P^T)^T(G+E)(Z_1+Y_2P)=0$.
Using (\ref{c2}), we obtain
$$
P(S_1^TS_1)=Y_2^TEZ_1-PZ_1^TEZ_1+Y_2^TEY_2P-PZ_1^TEY_2P,
$$
from which we obtain
$$
P=\epsilon \overline V_2^T(X_1^\dag X_2)^TV_1+{\cal O}(\epsilon^3).
$$
From (\ref{3.01})-(\ref{3.02}), $[\tilde Z_1\quad \tilde Y_2]$ has the following form
$$
[\tilde Z_1\quad \tilde Y_2]=\left[\begin{array}{cc}
V_1+{\cal O}(\epsilon^2)&-\epsilon B\overline V_2+{\cal O}(\epsilon^3)\\
\epsilon B^TV_1+{\cal O}(\epsilon^3)&\overline V_2+{\cal O}(\epsilon^2)\end{array}\right],
$$
and the representations of $G_\epsilon$ with respect to $\tilde Z_1, \tilde Y_2$ are given by $\tilde L_1=S_1^TS_1+{\cal O}(\epsilon^2),$  and
$$
\begin{array}{rl}
 \tilde L_2&=(E_{22}-PE_{12})(1+{\cal O}(\epsilon^2))=(Y_2^TEY_2-PZ_1^TEY_2)(1+{\cal O}(\epsilon^2))\\
&=\epsilon^2\overline V_2^T[ X_2^TX_2-(X_1^\dag X_2)^T(X_1^T X_2)]\overline V_2+{\cal O}(\epsilon^4)\\
 &=\epsilon^2\overline V_2^T\overline X_2^T\overline X_2\overline V_2+{\cal O}(\epsilon^4)=\epsilon^2 \overline S_2^T\overline S_2+{\cal O}(\epsilon^4).
 \end{array}
$$
Note that $G_\epsilon=X_\epsilon^TX_\epsilon$ is symmetric and has $\tilde Z_1$ and $\tilde Y_2$ as the bases of its right and left simple invariant subspaces such that
$\tilde Y_2^TG_\epsilon \tilde Z_1=0$, therefore $H:=[\tilde Z_1\quad  \tilde Y_2]$ satisfies
 $$
 H^TX_\epsilon^TX_\epsilon H=\left[\begin{array}{cc}S_1^TS_1+{\cal O}(\epsilon^2)& 0\\0 &\epsilon^2 \overline S_2^T\overline S_2+{\cal O}(\epsilon^4)\end{array}\right],
 $$
in which the orthonormal columns of $V_\epsilon:=[\tilde Z_1\quad \tilde Y_2]$ span the right singular subspace of $X_\epsilon$, with   diagonal entries of $S_1+{\cal O}(\epsilon^2), \epsilon \overline S_2+{\cal O}(\epsilon^3)$  as its singular values. It's obvious that
 $$
 X_\epsilon[\tilde Z_1\quad \tilde Y_2]=[X_1V_1+{\cal O}(\epsilon^2)\quad \epsilon \overline X_2 \overline V_2 ]+{\cal O}(\epsilon^3)=[U_1S_1+{\cal O}(\epsilon^2)\quad \epsilon \overline U_2 \overline S_2+{\cal O}(\epsilon^3)],
 $$
 from which we conclude that the left singular matrix $U_\epsilon$ of $X_\epsilon$ satisfies
 $$
U_\epsilon=\left[ U_1+{\cal O}(\epsilon^2) \quad \overline U_2+{\cal O}(\epsilon^2)\right].
 $$
The proof is then complete.\qed

{\it {\bf Theorem 3.3} For the multidimensional TLSE problem (\ref{1.2}),  with the notations in (\ref{2.03})--(\ref{2.6}), assume that $\overline V_{22}$ has full row rank, and the minimum Frobenius norm solution $X_t=-\overline V_{12}\overline V_{22}^\dag$. Denote \begin{equation}
L_\epsilon=W_{\epsilon}^{-1}L=\left[\begin{array}c \epsilon^{-1}C\\ A\end{array}\right],\qquad
H_\epsilon=W_{\epsilon}^{-1}H=\left[\begin{array}c \epsilon^{-1}D\\ B\end{array}\right],\label{3.3}
\end{equation}
where $W_\epsilon={\rm diag}(\epsilon I_{p}, I_q)$ with $\epsilon$ being a small positive parameter. Consider the multidimensional weighted TLS problem
\begin{equation}
\min_{\bar E, \bar f}\|[~\bar E\quad \bar f~]\|_F \qquad \mbox{subject to}\qquad (L_\epsilon+\bar E)X_\epsilon=H_\epsilon+\bar F,\label{3.4}
\end{equation}
then the minimum Frobenius norm solution $X_{t(\epsilon)}$ tends to $X_{t}$ as $\epsilon$ tends to zero.
}

{ {\it Proof.}}  To prove the close relation of TLSE solution to the WTLS solution, we need to investigate the right singular vectors of $\widetilde L_\epsilon=[L_\epsilon\quad H_\epsilon]$ corresponding to small singular values, in which $\widetilde L_\epsilon^T$
has the same left singular vectors as $[\widetilde C^T\quad \epsilon \widetilde A^T]$, and their  singular values are identical  up to multiplication by $\epsilon^{-1}$.

To apply Lemma 3.2, let  $\widetilde C^T=V_CS_CU_C^T$ be the skinny  SVD  of the full column-rank matrix $\widetilde C^T$, and the   SVD of $\widetilde A\widetilde Q_2$ be given by (\ref{2.4}).
It is obvious that
$$(I_{n+d}-\widetilde C^T\widetilde C^{\dag T})\widetilde A^T=\widetilde Q_2\widetilde Q_2^T\widetilde A^T=(\widetilde Q_2\widetilde V)\widetilde \Sigma^T\widetilde U^T=\overline V\widetilde\Sigma^T \widetilde U^T. $$
By Lemma 3.2, we know that the left and right singular matrices $\widetilde V_{\epsilon}$, $\widetilde U_{\epsilon}$ of $[\widetilde C^T\quad  \epsilon \widetilde A^T]$  satisfies
\begin{equation}
\begin{array}l
\widetilde V_{\epsilon}=\left[\underbrace{V_C+{\cal O}(\epsilon^2)}~~  \quad \underbrace{\overline V+{\cal O}(\epsilon^2)}\right],\qquad
\widetilde U_{\epsilon}=  \left[\underbrace{{\bf P}_\epsilon \Big(U_C+{\cal O}(\epsilon^2)}\Big)~~~~~ \underbrace{{\bf Q}_\epsilon\Big(\widetilde U+{\cal O}(\epsilon^2)\Big)}\right],\\
\qquad\qquad\quad {p\qquad ~~~~n+d-p}\qquad\qquad\qquad\qquad~~~~  p\qquad\quad\quad ~~~~~~n+d-p
\end{array}\label{3.5}
\end{equation}
where $\overline V=\widetilde Q_2\widetilde V_2$,
\begin{equation}
{\bf P }_\epsilon=\left[\begin{array}c I_p\\ \epsilon (\widetilde A\widetilde C^\dag)\end{array}\right],\qquad {\bf Q }_\epsilon=\left[\begin{array}c -\epsilon (\widetilde A\widetilde C^\dag)^T\\ I_q\end{array}\right],\label{3.05}
\end{equation}
and the corresponding  singular values are just   diagonal entries  of $S_C+{\cal O}(\epsilon^2)$,
$\epsilon\widetilde \Sigma+{\cal O}(\epsilon^3)$. Therefore the SVD of $\widetilde L_\epsilon$ is given by
$\widetilde L_\epsilon=\widetilde U_\epsilon \widetilde S_\epsilon\widetilde V_\epsilon^T$ for
\begin{equation}
\widetilde S_\epsilon={\rm diag}(\epsilon^{-1}S_C+{\cal O}(\epsilon), \widetilde \Sigma+{\cal O}(\epsilon^2)),\label{3.6}
\end{equation}
and  the smallest $n+d-p$  singular values of $\widetilde L_\epsilon$ can be approximated by
 $\widetilde\sigma_i+{\cal O }(\epsilon^2)$ for $i=1,\ldots, n+d-p$, and for sufficiently small $\epsilon$,
 $$
\widetilde\sigma_1+{\cal O }(\epsilon^2)\ge  \cdots\ge \widetilde\sigma_k+{\cal O }(\epsilon^2)>\widetilde\sigma_{k+1}+{\cal O }(\epsilon^2)\ge \ldots\ge \widetilde\sigma_{n+d-p}+{\cal O }(\epsilon^2),
 $$
 and the bottom  right $d\times (n+d-t)$ submatrix in $\widetilde V_\epsilon$ has  full row rank.
 Therefore the minimum Frobenius norm WTLS solution $X_{t(\epsilon)}$  to problem  (\ref{3.4}), in the limit, takes the form
$$
\lim_{\epsilon\rightarrow 0+} X_{t(\epsilon)}=\lim_{\epsilon\rightarrow 0+}\Big[-(\overline V_{12}+{\cal O}(\epsilon^2))(\overline V_{22}+{\cal O}(\epsilon^2))^\dag\Big]=-\overline V_{12}\overline V_{22}^\dag,
$$
 which is exactly $ X_{t}$. The proof of the theorem then follows.\qed\\

\section{Condition numbers of TLSE}

Condition numbers  measure the sensitivity of the solution to the original data in problems, and they play an
important role in numerical analysis.

To evaluate the condition number of the multidimensional TLSE problem, let $m=p+q$,
  $[\hat L\quad \hat H]=[L\quad H]+[\Delta L\quad \Delta H]$, where  the perturbation $[\Delta L\quad \Delta H]$ is sufficiently small.
In order to derive the first order perturbation estimate of  the TLSE solution, we
 define the mapping   $\phi: {\mathbb R}^{m(n+d)}\rightarrow {\mathbb R}^{nd}$ for the multidimensional TLSE problem  (\ref{1.2}):
 $$
\phi(c)= {\rm vec}(X_t),\quad c={\rm vec}([L\quad H]).
 $$
Define the absolute normwise,
relative normwise, mixed  and componentwise condition numbers of $X_{t}$   as follows
$$
\begin{array}l
\kappa^{\rm abs}(X_{t}, L, H)=\lim\limits_{\epsilon\rightarrow 0}\sup \left\{
{\displaystyle\|\Delta X_t\|_F\over\displaystyle \|[\Delta L\quad \Delta H]\|_F}: \|[\Delta L\quad \Delta H]\|_F\le \epsilon \|[ L\quad H]\|_F\right\},\\
\kappa^{\rm rel}(X_{t}, L, H)=\lim\limits_{\epsilon\rightarrow 0}\sup \left\{
{\displaystyle\|\Delta X_t\|_F\over\displaystyle \epsilon\|X_{t}\|_F}: \|[\Delta L\quad \Delta H]\|_F\le \epsilon \|[ L\quad H]\|_F \right\},\\
m(X_{t}, L, H)=\lim\limits_{\epsilon\rightarrow 0}\sup \left\{
{\displaystyle\|\Delta X_t\|_{\rm max}\over\displaystyle  \epsilon\|X_{t}\|_{\rm max}}: |\Delta L|\le \epsilon|L|,\quad |\Delta H]|\le \epsilon |H|\right\},\\
c(X_{t}, L, H)=\lim\limits_{\epsilon\rightarrow 0}\sup \left\{
{\displaystyle 1\over\displaystyle \epsilon}\left\|{\displaystyle\Delta X_t\over\displaystyle X_{t}}\right\|_{\rm max}: |\Delta L|\le \epsilon|L|,\quad |\Delta H]|\le \epsilon |H|\right\},
\end{array}
$$
where $|\cdot|$ denotes the componentwise absolute value, $Y\le Z$ means $y_{ij}\le z_{ij}$ for all $i,j$, and ${Y\over Z}$ is the  entry-wise division defined by ${Y\over Z}:=[{y_{ij}\over z_{ij}}]$ and  ${\xi\over 0}$
 is interpreted as zero if $\xi=0$ and infinity otherwise.

If  ${\rm vec}(X_t)=\phi(c)$ is  continuous and  Fr{\rm\'{e}}chet
differentiable  at the neighbourhood of the point $c$, according to the concept and formulae in \cite{ge, gk, rice}, the above condition numbers can be  formulated as follows:
$$
\begin{array}l
\kappa^{\rm abs}(X_{t}, L, H)={\|\phi'(c)\|_2},\qquad \kappa^{\rm rel}(X_{t}, L, H)=\frac{\displaystyle\|\phi'(c)\|_2\|c\|_2}{\displaystyle\|\phi(c)\|_2},\\
m(X_{t}, L, H)=\frac{\displaystyle\||\phi'(c)|\cdot |c|\|_{\infty}}{\displaystyle\|\phi(c)\|_{\infty}},\qquad
c(X_{t}, L, H)=\left\|{\displaystyle|\phi'(c)|\cdot|c|\over\displaystyle |\phi(c)|} \right\|_{\infty}.
\end{array}
$$

\subsection{Normwise condition number}

Notice that  $\phi'(c)$ is vital for above condition numbers, while a simple and  Fr{\rm\'{e}}chet
differentiable   expression of $\phi(c)$ is not easy to derive.
To get $\phi'(c)$, as did in \cite{liu2}, we start from the differentiability
of the weighted TLS solution $X_{t(\epsilon)}$ by defining  the mapping for the  multidimensional  WTLS  problem (\ref{3.3})-(\ref{3.4}):
$$
{\rm vec}(X_{t(\epsilon)})=\varphi(c_\epsilon),\qquad c_\epsilon={\rm vec}([L_\epsilon\quad  H_\epsilon]),
$$
and then get   the first order perturbation estimate ${\rm vec}(\Delta X_{t(\epsilon)})$ of WTLS solution based on the result in (\ref{2.3}), from which
  the first order perturbation estimate of the  TLSE solution is derived by taking the limit $\epsilon \rightarrow 0$.
 The idea of using limit technique to perform perturbation  and condition number analysis  of a problem was also used  by Wei and De Pierror \cite{pw,wp} for equality constrained least squares problem, and by Zheng and Yang \cite{zy} for mixed least squares-total least squares problem.

{\bf Lemma 4.1 }\cite{zmw} Let ${\small Q=\left[\begin{array}{cc} Q_{11}&Q_{12}\\Q_{21}&Q_{22}\end{array}\right]}$ be an $n$-by-$n$ orthogonal matrix with a 2-by-2 partitioning, then

(a) $Q_{11}$ has full column (row) rank if and only if $Q_{22}$ has full row (column) rank;

(b) $\|Q_{11}^\dag\|_2=\|Q_{22}^\dag\|_2,\quad Q_{11}^{\dag T}=Q_{11}-Q_{12}Q_{22}^\dag Q_{21},\quad Q_{11}^{\dag T}Q_{21}^T=-Q_{12}Q_{22}^\dag$.\\

 {\it{\bf Theorem 4.2} With the notation in (\ref{2.03})-(\ref{2.4}), let the skinny SVD of $\widetilde C$ be $\widetilde C= U_CS_CV_C^T$  and assume that the condition (\ref{2.5}) holds with the  partition $\overline V_{22}$ in  (\ref{2.6}) of full row rank. Denote
 $$
 \begin{array}l
{\bf P}=\left[\begin{array}c I_p\\ 0_{q\times p}\end{array}\right],\qquad {\bf Q}=\left[\begin{array}c -(\widetilde A\widetilde C^\dag)^T \\ I_q\end{array}\right],\\[10pt]
S_1=\left[\begin{array}{cc}
S_C&0\\0&\widetilde \Sigma_1\end{array}\right],\qquad
\widehat V_1=[V_C\quad \overline V_1 ]=\begin{array}l n\\d\end{array}\left[\begin{array} c\widehat V_{11}\\\hline \widehat V_{21}\end{array}\right].\\
\qquad\quad~~~~\quad~  \qquad \qquad\qquad\quad  p\quad~~ k\qquad\qquad~~~~ t \\
\end{array}
$$
Then for sufficiently small perturbation $\|[\Delta L\quad \Delta H]\|_F$, the first order perturbation estimate for the minimum  Frobenius norm TLSE solution $X_{t}=-\overline V_{12}\overline V_{22}^\dag$ takes the form
\begin{equation}
{\rm vec}(\Delta X_t)=K{\rm vec}([\Delta L\quad \Delta H])+{\cal O}(\|[\Delta L\quad \Delta H]\|_F^2),\label{4.4}
 \end{equation}
 where   $K=(H_1+H_2)G\widehat Z$ is exactly the Fr$\acute{e}$chet derivative $\phi'(c)$  and with $F_{ \overline V_{22}}=I- \overline V_{22}^\dag \overline  V_{22}$,
 \begin{equation}
\begin{array}l
H_{1}=\Big((\overline  V_{22}\overline V_{22}^T)^{-1}\widehat V_{21}\Big)\otimes (\overline  V_{12}F_{\overline V_{22}}),  \qquad H_{2}=\Big(\overline V_{22}^{\dag^T}\otimes \widehat  V_{11}^{\dag^T}\Big)\Pi_{(n+d-t,t)}, \\
G=(  S_{1}^2\otimes I_{n+d-t}-\left[\begin{array}{cc} 0_p&0\\ 0& I_k\end{array}\right]\otimes ( \widetilde  \Sigma_{2}^T \widetilde \Sigma_{2}))^{-1}\big[ I_t\otimes \widetilde \Sigma_{2}^T\qquad S_{1}\otimes I_{n+d-t}\big],\\
\widehat Z=\left[\begin{array}c  [0_{(n+d)\times p}\quad \overline V_1]^T\otimes({\bf Q} \widetilde U_{2})^T\\ \Pi_{(t,n+d-t)}
\Big(\overline V_2^T\otimes [{\bf P}U_C\quad {\bf Q}\widetilde U_1]^T\Big) \end{array}\right].
\end{array}\label{4.02}
\end{equation}}

{\it Proof.} Assume that the SVD of $\widetilde L_\epsilon=[L_\epsilon\quad H_\epsilon]=\widetilde U_\epsilon \widetilde S_\epsilon \widetilde V_\epsilon^T$ is given by  (\ref{3.5})--(\ref{3.6}), whose factors have partitions as
\begin{equation}
\begin{array}l
 \widetilde U_\epsilon=\Big[ \widetilde U_{1(\epsilon)}  ~|~ \widetilde U_{2(\epsilon)}\Big]=\Big[{\bf P}_\epsilon U_C\quad {\bf Q}_\epsilon\widetilde U_1~~\Big|~~{\bf Q}_\epsilon\widetilde U_2\Big]+{\cal O}(\epsilon^2),\\
 \qquad \qquad{t\quad~ n+d-t\quad ~ p~~~~~~~~~~~k~~~~~~~~~~n+d-t}\\
\widetilde V_\epsilon=\begin{array}l n\\ d\end{array}\left[\begin{array}{c|c}
\widetilde V_{11(\epsilon)}~&~ \widetilde V_{12(\epsilon)}\\
\widetilde V_{21(\epsilon)}~&~ \widetilde V_{22(\epsilon)}\end{array}
\right]=\left[\underbrace {V_C\quad \overline V_1}~~ \Big|~~ \overline V_2\right]+{\cal O}(\epsilon^2),\\
 \qquad \qquad{\quad t\quad~~~~ n+d-t} \qquad \qquad ~~{t}   \qquad~ n+d-t\\
 \widetilde S_{1(\epsilon)}={\rm diag}(\epsilon^{-1}S_C+{\cal O}(\epsilon),~ \widetilde \Sigma_1+{\cal O}(\epsilon^2)),\quad  \widetilde S_{2(\epsilon)}=\widetilde \Sigma_2+{\cal O}(\epsilon^2),
 \end{array}\label{4.002}
\end{equation}
in which ${\bf P}_\epsilon, {\bf Q}_\epsilon$ are given by (\ref{3.05}).

By applying the result in (\ref{2.3}) for the WTLS problem (\ref{3.3})--(\ref{3.4}),  the first order perturbation estimate of the WTLS solution $X_{t(\epsilon)}$   satisfies
\begin{equation}
\begin{array}{rl}
{\rm vec}(\Delta X_{t(\epsilon)})&=\varphi'(c_{\epsilon}){\rm vec}([\Delta L_\epsilon\quad \Delta H_\epsilon])+{\cal O}(\|\Delta L_\epsilon\|_F^2+ \|\Delta H_\epsilon\|_F^2)\\
&=(H_{1(\epsilon)}+H_{2(\epsilon)})D_\epsilon Z_\epsilon{\rm vec}([\Delta L_\epsilon\quad  \Delta H_\epsilon])+{\cal O}(\|\Delta L_\epsilon\|_F^2+ \|\Delta H_\epsilon\|_F^2),
\end{array}\label{4.3}
\end{equation}
where
\begin{equation}
\begin{array}l
H_{1(\epsilon)}=\Big((\widetilde  V_{22(\epsilon)}\widetilde V_{22(\epsilon)}^T)^{-1}\widetilde V_{21(\epsilon)}\Big)\otimes (\widetilde  V_{12(\epsilon)}F_{\widetilde V_{22(\epsilon)}}),\\
H_{2(\epsilon)}=\Big(\widetilde V_{22(\epsilon)}^{\dag^T}\otimes \widetilde  V_{11(\epsilon)}^{\dag^T}\Big)\Pi_{(n+d-t,t)},\\
D_\epsilon=(\widetilde S_{1(\epsilon)}^2\otimes I_{n+d-t}-I_t\otimes ( \widetilde S_{2(\epsilon)}^T \widetilde S_{2(\epsilon)}))^{-1}\big[ I_t\otimes \widetilde S_{2(\epsilon)}^T\qquad\widetilde S_{1(\epsilon)}\otimes I_{n+d-t}\big],\\
Z_\epsilon=\left[\begin{array}c\widetilde V_{1(\epsilon)}^T\otimes\widetilde U_{2(\epsilon)}^T\\ \Pi_{(t,n+d-t)}(\widetilde V_{2(\epsilon)}^T\otimes \widetilde U_{1(\epsilon)}^T)\end{array}\right],
\end{array}\label{4.01}
\end{equation}
with $F_{ \widetilde V_{22(\epsilon)}}=I- \widetilde V_{22(\epsilon)}^\dag \widetilde  V_{22(\epsilon)}$.

In (\ref{4.3}), denote $Y_\epsilon=D_\epsilon Z_\epsilon{\rm vec}([\Delta L_\epsilon\quad  \Delta H_\epsilon])$. Notice that $[\Delta L_\epsilon\quad  \Delta H_\epsilon]=W_\epsilon^{-1}[\Delta L\quad  \Delta H]$ for $W_\epsilon={\rm diag}(\epsilon I_p, I_q)$, therefore
\begin{equation}
Y_\epsilon=D_\epsilon\left[\begin{array}{c}
{\rm vec}\Big(\widetilde U_{2(\epsilon)}^T[\Delta L_\epsilon\quad  \Delta H_\epsilon]\widetilde V_{1(\epsilon)}\Big)\\
{\rm vec}\Big(\widetilde V_{2(\epsilon)}^T[\Delta L_\epsilon\quad  \Delta H_\epsilon]^T\widetilde U_{1(\epsilon)}\Big)
\end{array}\right]=D_\epsilon\left[\begin{array}{c}
{\rm vec}\Big(\widetilde U_{2(\epsilon)}^TW_\epsilon^{-1}[\Delta L\quad  \Delta H]\widetilde V_{1(\epsilon)}\Big)\\
{\rm vec}\Big(\widetilde V_{2(\epsilon)}^T[\Delta L\quad  \Delta H]^T(W_\epsilon^{-1}\widetilde U_{1(\epsilon)})\Big)
\end{array}\right].\label{4.2}
\end{equation}
Set $\widehat W_\epsilon={\rm diag}(\epsilon I_p, I_k)$ and
$\widehat S_{1(\epsilon)}=\widehat W_\epsilon\widetilde S_{1(\epsilon)}$, then  $\widehat S_{1(\epsilon)}=S_1+{\cal O}(\epsilon^2)),$
and   (\ref{4.2}) becomes
$$
\begin{array}{rl}
Y_\epsilon&= D_\epsilon\left[\begin{array}{c}
{\rm vec}\Big(\widetilde U_{2(\epsilon)}^TW_\epsilon^{-1}[\Delta L\quad  \Delta H]\widetilde V_{1(\epsilon)}\Big)\\
{\rm vec}\Big(\widetilde V_{2(\epsilon)}^T[\Delta L\quad  \Delta H]^T(W_\epsilon^{-1}\widetilde U_{1(\epsilon)})\Big)
\end{array}\right]= G_\epsilon\widehat Z_\epsilon{\rm vec}([\Delta L\quad  \Delta H]),
\end{array}
$$
where
$$
\begin{array}l
 D_\epsilon=(\widehat S_{1(\epsilon)}^{~2}\otimes I_{n+d-t}-\widehat W_\epsilon^2\otimes ( \widetilde S_{2(\epsilon)}^T \widetilde S_{2(\epsilon)}))^{-1}\big[ \widehat W_\epsilon^2\otimes \widetilde S_{2(\epsilon)}^T\qquad(\widehat W_\epsilon\widehat S_{1(\epsilon)})\otimes I_{n+d-t}\big],\\
G_\epsilon=(\widehat S_{1(\epsilon)}^{~2}\otimes I_{n+d-t}-\widehat W_\epsilon^2\otimes ( \widetilde S_{2(\epsilon)}^T \widetilde S_{2(\epsilon)}))^{-1}\big[ I_t\otimes \widehat S_{2(\epsilon)}^T\qquad \widehat S_{1(\epsilon)}\otimes I_{n+d-t}\big],\\
\widehat Z_\epsilon=\left[\begin{array}c
(\widetilde V_{1(\epsilon)}\widehat W_\epsilon^2)^T\otimes (\widetilde U_{2(\epsilon)}^TW_\epsilon^{-1})\\
\Pi_{(t,n+d-t)}\Big( \widetilde V_{2(\epsilon)}^T\otimes (W_\epsilon^{-1}\widetilde U_{1(\epsilon)}\widehat W_\epsilon)^T\Big) \end{array}\right].
\end{array}
$$
By the expressions in (\ref{4.002}), and    taking the limit $\epsilon\rightarrow 0$ for $H_{1(\epsilon)}, H_{2(\epsilon)}, G_{\epsilon}$ and $\widehat Z_{\epsilon}$ in  (\ref{4.01}), we    obtain
the corresponding limit matrices $H_1, H_2, G, \widehat Z$ as  (\ref{4.02}), and $K=(H_1+H_2)G\widehat Z$ is exactly the  Fr{\rm\'{e}}chet derivative $\phi'(c)$.\qed\\

{\it{\bf Theorem 4.3} With the notation in Theorem 4.2, the absolute and relative condition numbers of the minimum Frobenius norm TLSE solution $X_t$ are given by
$$
\kappa^{\rm abs}(X_t, L, H)=\|(H_1+H_2)G\overline Z\|_2,\qquad \kappa^{\rm rel}(X_t, L, H)=\|(H_1+H_2)G\overline Z\|_2{\|[L\quad  H]\|_F\over \|X_t\|_F},
$$
where
$$
\overline Z={\rm diag}\Big(\begin{array}{cc}\begin{array}{cc}
\left[\begin{array}{cc} 0_p&0\\ 0& I_k\end{array}\right]\otimes (\widetilde U_2^T{\bf Q}^T)\end{array},\quad
\left[\begin{array}{cc} I_p&~0\\-\widetilde U_1^T(\widetilde A\widetilde C^\dag)U_C&~I_k\end{array}\right]\otimes I_{n+d-t}\end{array}\Big).
$$
In particular, when $k=n-p$, the term $H_1$ diminishes to zero, and $H_{2}=\big(\overline V_{22}^{-T}\otimes \widehat  V_{11}^{-T}\big)\Pi_{(d,n)}.$
}

{\it Proof.} By the condition number formulae,  the absolute and relative condition numbers of the solution $X_t$  are given by
$$
\begin{array}l
\kappa^{\rm abs}(X_t, L, H)=\|\phi'(c)\|_2=\|K\|_2,\quad \\
\kappa^{\rm rel}(X_t, L, H)={\displaystyle\|\phi'(c)\|_2\|c\|_2\over\displaystyle \|X_t\|_{F}}={\displaystyle \|K\|_2\|[L\quad H]\|_F\over \displaystyle \|X_t\|_{F }},
\end{array}
$$
in which
$$
\|K\|_2=\|KK^T\|_2^{1/2}=\|(H_1+H_2)G\widehat Z\widehat Z^TG^T(H_1+H_2)\|^{1/2},
$$
for
$$
\widehat Z\widehat Z^T=\left[\begin{array}{cc}\begin{array}{cc}
\left[\begin{array}{cc} 0_p&0\\ 0& I_k\end{array}\right]\otimes (\widetilde U_2^T{\bf Q}^T{\bf Q}\widetilde U_2)\end{array}&0\\
0&\Pi_{(t,n+d-t)}\Big(I_{n+d-t}\otimes M\Big)\Pi_{(t,n+d-t)}^T\end{array}\right]=\breve Z\breve Z^T,
$$
and
$$
\begin{array}l
M=\left[\begin{array}{cc} I_p&U_C^T{\bf P}^T{\bf Q}\widetilde U_1\\ \widetilde U_1^T{\bf Q}^T{\bf P}U_C&\widetilde U_1^T{\bf Q}^T{\bf Q}\widetilde U_1\end{array}\right]=
\left[\begin{array}{cc} I_p&~0\\-\widetilde U_1^T(\widetilde A\widetilde C^\dag)U_C&~I_k\end{array}\right]\left[\begin{array}{cc} I_p&~0\\-\widetilde U_1^T(\widetilde A\widetilde C^\dag)U_C&~I_k\end{array}\right]^T,\\
\breve Z={\rm diag}\Big(\begin{array}{cc}\begin{array}{cc}
\left[\begin{array}{cc} 0_p&0\\ 0& I_k\end{array}\right]\otimes (\widetilde U_2^T{\bf Q}^T)\end{array},\quad
\Pi_{(t,n+d-t)}\Big(I_{n+d-t}\otimes \left[\begin{array}{cc} I_p&~0\\-\widetilde U_1^T(\widetilde A\widetilde C^\dag)U_C&~I_k\end{array}\right]\Big)\end{array}\Big)\\
\quad={\rm diag}\Big(\begin{array}{cc}\begin{array}{cc}
\left[\begin{array}{cc} 0_p&0\\ 0& I_k\end{array}\right]\otimes (\widetilde U_2^T{\bf Q}^T)\end{array},\quad
\Big( {\left[\begin{array}{cc} I_p&~0\\-\widetilde U_1^T(\widetilde A\widetilde C^\dag)U_C&~I_k\end{array}\right]\otimes I_{n+d-t}}\Big)\Pi_{(n+d-t,t)}\end{array}\Big).
\end{array}
$$
Notice that  $\Pi_{(n+d-t,t)}$  is an orthogonal matrix, then
$$
\|K\|_2=\|(H_1+H_2)G\breve Z\|_2=\|(H_1+H_2)G\overline Z\|_2.
$$

When $k=n-p$, $\overline V_{22}$ is nonsingular and $X_{n}=-\overline V_{12}\overline V_{22}^{-1}$. Moreover,
 note that the SVDs of $\widetilde A\widetilde Q_2\widetilde Q_2^T=\widetilde U\widetilde \Sigma \overline  V^T$ and $\widetilde C=U_CS_CV_C^T$
imply ${\cal R}(V_C)={\cal R}(\widetilde C^T)={\cal R}(\widetilde Q_1)$ and ${\cal R}(\overline V)\subseteq {\cal R}(\widetilde Q_2)$,  therefore $V_C^T\overline V=0$ and
\begin{equation}
\begin{array}l
\breve V:=[V_C\quad \overline V]=\left[\begin{array}{cc}
\widehat V_{11}&~~\overline V_{12}\\
\widehat V_{21}&~~\overline V_{22}\end{array}\right]
\end{array}\label{4.001}
\end{equation}
is an $(n+d)\times (n+d)$ orthogonal matrix. According to Lemma 4.1(a),
 $\widehat V_{11}$ is nonsingular. It follows that
 $H_1=0$, $H_{2}=\big(\overline V_{22}^{-T}\otimes \widehat  V_{11}^{-T}\big)\Pi_{(d,n)}.$  This completes the proof.\qed

{\it {\bf Theorem 4.4} Let
$$
 \rho_{AC}^{(1)}=1+\|\widetilde C\|_2+\|\widetilde A\widetilde C^\dag \widetilde C\|_2,\qquad \rho_{AC}^{(2)}=1+\|\widetilde C^\dag\|_2+\|\widetilde A\widetilde C^\dag\|_2, \qquad  \eta_k^\sigma= \max\{1,{\displaystyle \sqrt{\widetilde\sigma_k^2+\widetilde\sigma_{k+1}^2}\over\displaystyle \widetilde\sigma_k^2-\widetilde\sigma_{k+1}^2}\},
$$
then for the absolute normwise condition number,  we have
$$
\kappa^{\rm abs}( X, L, H)\le (1+\|X_t\|_2^2)\rho_{AC}^{(2)}\eta_k^{\sigma}.
$$
 In particular, when $k=n-p$, it has the  bounds as
$$
\begin{array}{rl}
{\displaystyle \eta_k^\sigma\over \displaystyle \|\widehat V_{11}\|_2\|\overline V_{22}\|_2\rho_{AC}^{(1)}}
&\le \kappa^{\rm abs}( X, L, H)\le (1+\|X_n\|_2^2)\rho_{AC}^{(2)}\eta_k^{\sigma}.
\end{array}
$$
}

{\it Proof.} Let $ W_0=\left[\begin{array}{cc} 0_p&0\\ 0& I_k\end{array}\right]$, it follows that $\overline Z=\Gamma \ddot{Z}$ for $\Gamma={\rm diag}(\Gamma_1,\Gamma_2)\otimes I_{n+d-t}$, and
$$
\Gamma_1={\rm diag}(W_0,~ I_{p+k})\otimes I_{n+d-t},\qquad \Gamma_2={\rm diag}(S_C, I_k)\otimes I_{n+d-t},
$$
and
$$
\ddot{ Z}={\rm diag}\Big(\begin{array}{cc}\begin{array}{cc}
W_0\otimes (\widetilde U_2^T{\bf Q}^T)\end{array},\quad
\left[\begin{array}{cc} S_C^{-1}&~0\\-\widetilde U_1^T(\widetilde A\widetilde C^\dag)U_C&~I_k\end{array}\right]\otimes I_{n+d-t}\end{array}\Big)=:{\rm diag}({\ddot Z}_{11},
{\ddot Z}_{22}).
$$
Therefore
\begin{equation}
\kappa^{\rm abs}(X_t, L, H)\le \|H_1+H_2\|_2\|\overline G\|_2\|\ddot{ Z}\|_2,\label{4a}
\end{equation}
where $\overline G=G\Gamma$, and
$$
\overline G\,\overline G^T=(\left[\begin{array}{cc} S_C^4&0\\ 0& \widetilde \Sigma_1^T\widetilde \Sigma_1\end{array}\right]\otimes I_{n+d-t}+W_0\otimes (\widetilde \Sigma_2^T\widetilde \Sigma_2) )(S_{1}^2\otimes I_{n+d-t}-W_0\otimes ( \widetilde  \Sigma_{2}^T \widetilde \Sigma_{2}))^{-2},
$$
consists of $(k+1)$  diagonal block $D^{(i)}$ for $i=0,1,\cdots,k$ satisfying
$$
D^{(0)}=I_{p(n+d-t)},\qquad  D^{(i)}={\rm diag}\Big({\displaystyle \widetilde\sigma_i^2+\widetilde\sigma_{k+j}^2\over(\displaystyle \widetilde\sigma_i^2-\widetilde\sigma_{k+j}^2)^2}\Big),
\quad 1\le i\le k, 1\le j\le n+d-t.
$$
Note that    ${\sigma^2+\eta^2\over
(\sigma^2-\eta^2)^2}$ is an increasing function  of $\eta$ and a decreasing function of $\sigma$ for $\sigma>\eta>0$, therefore
\begin{equation}
\|\overline G\|_2=\|\overline G\,\overline G^T\|_2^{1/2}=  \eta_k^\sigma.\label{4b}
\end{equation}

For the upper bound of $\|H_1+H_2\|_2$, note that $\breve V$ in (\ref{4.001})
is an   orthogonal matrix and   $X_t=-\overline V_{12}\overline V_{22}^\dag$, then by applying the CS decomposition (see \cite[Theorem 2.6.3]{gv2}) and a similar technique in \cite[Theorem 3.6]{zmw}, we have
\begin{equation}
\|H_1+H_2\|_2\le \|\overline V_{22}^\dag \|_2^2=1+\|X_t\|_2^2.\label{4c}
\end{equation}

For the norm of ${\ddot Z}$ and ${\ddot Z}^\dag$,  note that
 $$
 \|{\ddot Z}_{11}\|_2\le \|{\bf Q}\|_2= 1+\|\widetilde A\widetilde C^\dag\|_2,\quad
\|{\ddot Z}_{11}^\dag\|_2= (\sigma_{\rm min}({\bf Q}\widetilde U_2))^{-1}=
 (\sigma_{\rm min}(I+\widetilde U_2^T\widehat C^T\widehat C\widetilde U_2))^{-1/2}\le   1$$
  for $\widehat C=(\widetilde A\widetilde C^\dag)^T$. Moreover, with $(\widetilde A\widetilde C^\dag)U_CS_C=\widetilde AV_C$ and $V_CV_C^T= \widetilde C^\dag\widetilde C$,
$$
\begin{array}l
\|\ddot Z_{22}\|_2=\left\|\left[\begin{array}{cc} S_C^{-1}&~0\\-\widetilde U_1^T(\widetilde A\widetilde C^\dag)U_C^T&~I_k\end{array}\right]\right\|_2\le 1+\|\widetilde C^\dag\|_2+\|\widetilde A\widetilde C^\dag\|_2=\rho_{AC}^{(2)},\\
\|\ddot Z_{22}^{-1}\|_2=\left\|\left[\begin{array}{cc} S_CV_C^T&~0\\ \widetilde U_1^T\widetilde AV_CV_C^T&~I_k\end{array}\right]\right\|_2\le
1+\|\widetilde C\|_2+\|\widetilde A\widetilde C^\dag \widetilde C\|_2=\rho_{AC}^{(1)}.
\end{array}
$$
Therefore
\begin{equation}
\|\ddot Z\|_2\le \rho_{AC}^{(2)},
\qquad \|\ddot Z^\dag \|_2\le \rho_{AC}^{(1)}.
\label{4d}
\end{equation}
Substituting (\ref{4b})--({\ref{4d}}) into (\ref{4a}),
the upper bound for $\kappa^{\rm abs}(X_t, L, H)$  then follows.

When $k=n-p$,  the upper bound of  $\kappa^{\rm abs}(X_n,L,H)$ is obvious and for the lower bound,
note that
$$
\begin{array}{rl}
\kappa^{\rm abs}(X,L,H)&\ge  \sigma_{\rm min}(H_2)\|\overline G\ddot Z\|_2 \ge \sigma_{\rm min}(H_2)\sigma_{\rm min}(\ddot Z)\|\overline G\|_2\\
&=
{\displaystyle \eta_k^\sigma\over \displaystyle \|H_2^{-1}\|_2\|\ddot Z^{\dag}\|_2}\ge {\displaystyle \eta_k^\sigma\over \displaystyle \|\widehat V_{11}\|_2\|\overline V_{22}\|_2\rho_{AC}^{(1)}}.
\end{array}
$$
 The assertion then follows. \qed\\

 The upper bounds and lower bounds in Theorem 4.4 reduce to the ones in (\ref{bd1})-(\ref{bd2}) for the TLS problem, when $\widetilde C=0$ and $\eta_k^\sigma>1$.
Moreover, note that  $\|\overline V_{22}^\dag\|_2^2=1+\|X_t\|_2^2=1/\sigma_{\rm min}^2(\overline V_{22})$, it follows from Theorem 4.4  that the multidimensional TLSE problem might be ill-conditioned,  when
$\widetilde C$ is ill conditioned, or the gap between $\widetilde\sigma_k$ and $\widetilde\sigma_{k+1}$ or $\sigma_{\rm min}(\overline V_{22})$ is small.

When $d=1$, the absolute condition number $\kappa^{\rm abs}(X_t, L, H)$ has a compact form as follows.

{\bf Corollary 4.5} With   notations in Theorems  4.2 and 4.4, let $d=1$, $t=k+p$ and $\widetilde\sigma_k>\widetilde\sigma_{k+1}$ such that $\overline V_{22}\not=0$ holds for the approximate system $Ax\approx b$ subject to $Cx=\bar d$, then we have
$$
\begin{array}{rl}
\kappa^{\rm abs}(x_t, L, H)&={\displaystyle 1\over\displaystyle \|V_{22}\|_2^2}\|(\widehat V_{21}\otimes (\overline V_{12}+x_t\overline V_{22})+(\widehat V_{11}+x_t\widehat V_{21})\otimes \overline V_{22})G\overline Z\|_2\\
&\le (1+\|x_t\|_2^2)\rho_{AC}^{(2)}\eta_{k}^\sigma,
\end{array}
$$
where $\widetilde A=[A\quad b], \widetilde C=[C\quad \bar d]$.
In particular, if $k=n-p$, then
$$
\begin{array}{rl}
{\displaystyle \eta_{n-p}^\sigma\over \displaystyle \|\widehat V_{11}\|_2\|\overline V_{22}\|_2 \rho_{AC}^{(1)}}&\le \kappa^{\rm abs}(x_n, L, H)={\displaystyle 1\over\displaystyle\|V_{22}\|_2^2}\|(\widehat V_{11}+x_t\widehat V_{21})\otimes \overline V_{22})G\overline Z\|_2\\
&\le
(1+\|x_n\|_2^2)\rho_{AC}^{(2)} \eta_{n-p}^\sigma.
\end{array}
$$

{\it Proof.} Note that when $d=1$,  $\overline V_{22}\in {\mathbb R}^{1\times (n+1-t)}$ and $\overline V_{22}^\dag= {\overline V_{22}^T\over \|\overline V_{22}\|_2^2}$, $x_t=-\overline V_{12}\overline V_{22}^\dag=-\overline V_{12}\overline V_{22}^T/\|\overline V_{22}\|_2^2$.

$$
\begin{array}l
H_{1}=\Big((\overline  V_{22}\overline V_{22}^T)^{-1}\widehat V_{21}\Big)\otimes (\overline  V_{12}F_{\overline V_{22}})={1\over\|V_{22}\|_2^2}\Big(\widehat V_{21}\otimes (\overline V_{12}+x_k\overline V_{22})\Big), \\
 H_{2}=\Big(\overline V_{22}^{\dag^T}\otimes \widehat  V_{11}^{\dag^T}\Big)\Pi_{(n+1-t,t)}=\Pi_{(1,n)}(\widehat  V_{11}^{\dag^T}\otimes \overline V_{22}^{\dag^T} )={1\over\|\overline V_{22}\|_2^2}(\widehat  V_{11}^{\dag^T}\otimes \overline V_{22} ),
 \end{array}
$$
where $\widehat  V_{11}$ is the upper-left $n\times t$ submatrix of the orthogonal  matrix $\breve V$ given in (\ref{4.001}).
By Lemma 4.1(b), we have
\begin{equation}
\widehat V_{11}^{\dag^T}=\widehat V_{11}+x_t\widehat V_{21}.\label{4.00}
\end{equation}
 The formula for $\kappa^{\rm abs}$ then follows, in which  $H_1=0$ when $k=n-p$. The proof then follows.\qed

\subsection{Mixed and componentwise condition numbers}

For the mixed and componentwise condition numbers, we have the following results.

{\it {\bf Theorem 4.6} With the notation in Theorem 4.2,    then we have mixed and componentwise condition formulae of $X_t$ as follows:
\begin{eqnarray}
m(X_t, L, H)&=&{\displaystyle \left\||MN|{\rm vec}([|L|\quad |H|])\right\|_\infty\over\displaystyle  \|X_t\|_{\rm max}},\label{4.5}\\
c(X_t, L, H)&=&\left\|{\displaystyle  |MN|{\rm vec}([|L|\quad |H|])\over\displaystyle  {\rm vec}(|X_t|)}\right\|_{\infty},\label{4.6}
\end{eqnarray}
where $M=(H_1+H_2)D^{-1}$, $N=N_1+N_2$  for
$$
\begin{array}l
D= S_{1}^2\otimes I_{n+d-t}-\left[\begin{array}{cc} 0_p&0\\ 0& I_k\end{array}\right]\otimes ( \widetilde  \Sigma_{2}^T \widetilde \Sigma_{2}),\\
N_1=[0_{(n+d)\times p}\quad \overline V_1]^T\otimes ({\bf Q}\widetilde U_2\widetilde \Sigma_2)^T,\quad N_2=\Pi_{(t,n+d-t)} \Big(\overline V_2^T\otimes [{\bf P}U_CS_C\quad {\bf Q}\widetilde U_1\widetilde \Sigma_1 ]^T\Big).
\end{array}
$$

Moreover, they have  compact upper bounds as
\begin{eqnarray}
m^u(X_t, L, H)&=&{\displaystyle \left\||\widehat V_{11}^{\dag T}| Y^T|\overline V_{22}^\dag|+|\overline V_{12}F_{\overline V_{22}}|Y|\widehat V_{21}^T(\overline V_{22}\overline V_{22}^T)^{-1}|\right\|_{\rm max}\over\displaystyle  \|X_t\|_{\rm max}},\nonumber\\
c^u(X_t, L, H)&=& \left\|{\displaystyle|\widehat V_{11}^{\dag T}|Y^T|\overline V_{22}^\dag|+|\overline V_{12}F_{\overline V_{22}}|Y|\widehat V_{21}^T(\overline V_{22}\overline V_{22}^T)^{-1}|\over\displaystyle  X_t}\right\|_{\rm max},\nonumber
\end{eqnarray}
where $Y$ satisfies $D{\rm vec}(Y)={\rm vec}(\Upsilon)$ for
\begin{equation}
\Upsilon=|{\bf Q}\widetilde U_2\widetilde \Sigma_2|^T[|L|\quad |H|][0_{(n+d)\times p}\quad |\overline V_1|] +|\overline V_2^T|[|L|\quad |H|]^T|[{\bf P}U_CS_C\quad {\bf Q}\widetilde U_1\widetilde \Sigma_1]|,\label{4.9}
\end{equation}
from which the   $i$-th column of $y_i=Ye_i$ takes the form
\begin{equation}
y_i=(s_i^2I_{n+d-t}-\tau_i\widetilde\Sigma_2^T\widetilde\Sigma_2)^{-1}\Upsilon e_i.\label{4.004}
\end{equation}
Here $s_i$ is the $i$-th diagonal element of $S_1$ and  $\tau_i=1$ for $i>p$ and zero otherwise.

}

{\it Proof.}  By Theorem 4.2 and   the concept of condition numbers, the mixed and componentwise condition numbers of $X_t$ can be formulated
$$
\begin{array}{rl}
&m(X_t, L, H)={\displaystyle \||\phi'(c)|\cdot|c|\|_\infty\over \displaystyle \|\phi(c)\|_\infty}=
{\displaystyle \left\||(H_1+H_2)G\widehat Z|\cdot{\rm vec}([|L|\quad |H|)\right\|_\infty\over \displaystyle \|X_t\|_{\rm max}}\\
&={\displaystyle \left\|\left|M\Big(N_1+(S_{1}\otimes I_{n+d-t})\Pi_{(t,n+d-t)}
\big(\overline V_2^T\otimes [{\bf P}U_C\quad {\bf Q}\widetilde U_1]^T\big)\right|\Big){\rm vec}([|L|\quad |H|)\right\|_\infty\over \displaystyle \|X_t\|_{\rm max}}\\
&={\displaystyle \left\|\left|M\Big(N_1+\Pi_{(t,n+d-t)}(I_{n+d-t}\otimes S_{1})
\big(\overline V_2^T\otimes [{\bf P}U_C\quad {\bf Q}\widetilde U_1]^T\big)\right|\Big){\rm vec}([|L|\quad |H|)\right\|_\infty\over \displaystyle \|X_t\|_{\rm max}}\\
&={\displaystyle \left\||MN|{\rm vec}([|L|\quad |H|])\right\|_\infty\over\displaystyle  \|X_t\|_{\rm max}},
\end{array}
$$
which is bounded by
$$
\begin{array}{rl}
&|MN|{\rm vec}([|L|\quad |H|])\le (|H_1|+|H_2|) D^{-1}(|N_1|+|N_2|){\rm vec}([|L|\quad |H|])\\
&\le (|H_1|+|H_2|)D^{-1}{\rm vec}(\Upsilon)\\[5pt]
 &=\Big(|(\overline V_{22}^{\dag T}\otimes \widehat V_{11}^{\dag T})\Pi_{(n+d-t,t)}|+|\big((\overline V_{22}\overline V_{22}^T)^{-1}\widehat V_{21}\big) \otimes (\overline V_{12}F_{\overline V_{22}})\Big) {\rm vec}(Y)\\[5pt]
 &\le {\rm vec}\Big(|\widehat V_{11}^{\dag T}|Y^T|\overline V_{22}^\dag|+|\overline V_{12}F_{\overline V_{22}}|Y|\widehat V_{21}^T(\overline V_{22}\overline V_{22}^T)^{-1}|\Big),
\end{array}
$$
where $D{\rm vec}(Y)={\rm vec}(\Upsilon)$ gives
$$
YS_1^2-\widetilde\Sigma_2^T\widetilde\Sigma_2Y\left[\begin{array}{cc} 0_p&0\\ 0& I_k\end{array}\right]=\Upsilon,
$$
from which  the upper bound $m^u(X_t, L, H)$ follows with the $i$-th column of $Y$ satisfying (\ref{4.004}).

For the componentwise condition number, it follows that
$$
c(X_t, L, H)=\left\|{\displaystyle |\phi'(c)|\cdot|c|\over \displaystyle \phi(c)}\right\|_\infty=\left\|{\displaystyle  |MN|{\rm vec}([|L|\quad |H|])\over\displaystyle  {\rm vec}(|X_t|)}\right\|_{\infty},
$$
and the upper bound   $c^u(X_t, L, H)$  follows obviously.  We then finish the proof of the theorem.\qed

The following result is straightforward from Theorem 4.6.

{\it {\bf Corollary 4.7} With notations in  Theorem 4.2, if for $k=n-p$,  $\widetilde \sigma_{k}>\widetilde \sigma_{k+1}$  and $\overline V_{22}$ is nonzero.
then  the mixed and componentwise condition numbers for the TLSE solution $X_n=-\overline V_{12}\overline V_{22}$ satisfies
\begin{eqnarray}
m(X_{n}, L, H)&=&{\displaystyle \left\||MN|{\rm vec}([|L|\quad |H|])\right\|_\infty\over\displaystyle  \|X_n\|_{\rm max}},\nonumber \\
c(X_{n}, L, H)&=&\left\|{\displaystyle  |MN|{\rm vec}([|L|\quad |H|])\over\displaystyle  {\rm vec}(|X_n|)}\right\|_{\infty}, \nonumber
\end{eqnarray}
where
\begin{eqnarray}
M&=&\big(\overline V_{22}^{-T}\otimes \widehat  V_{11}^{-T}\big)\Pi_{(d,n)}(  S_{1}^2\otimes I_{d}-\left[\begin{array}{cc} 0_p&0\\ 0& I_{n-p}\end{array}\right]\otimes ( \widetilde  \Sigma_{2}^T \widetilde \Sigma_{2}))^{-1},\label{4.7}\\
N&=&[0_{(n+d)\times p}\quad \overline V_1]^T\otimes ({\bf Q}\widetilde U_2\widetilde \Sigma_2)^T+
\Pi_{(n,d)} \Big(\overline V_2^T\otimes [{\bf P}U_CS_C\quad {\bf Q}\widetilde U_1\widetilde \Sigma_1 ]^T\Big).\label{4.8}
\end{eqnarray}
Moreover, they have   upper bounds as
\begin{equation}
m^u(X_n, L, H)={\displaystyle \left\||\overline V_{11}^{-T}| Y^T|\overline V_{22}^{-1}|\right\|_{\rm max}\over\displaystyle  \|X_n\|_{\rm max}},\qquad
c^u(X_n, L, H)= \left\|{\displaystyle|\overline V_{11}^{-T}|Y^T|\overline V_{22}^{-1}|\over\displaystyle  X_n}\right\|_{\rm max}, \nonumber
\end{equation}
where the $i$-th column of $y_i=Ye_i$ satisfies
\begin{equation}
y_i=(s_i^2I_{d}-\tau_i\widetilde\Sigma_2^T\widetilde\Sigma_2)^{-1}\Upsilon e_i,\nonumber
\end{equation}
in which $s_i$ is the $i$-th diagonal element of $S_1$ and  $\tau_i=1$ for $i>p$ and zero otherwise, $\Upsilon$ is given by  (\ref{4.9}) with $t=n$. \\

}

{\it {\bf Theorem 4.8} For the single-dimensional TLSE problem ($d=1$), assume that for $k=n-p$,
$\widetilde \sigma_{k}>\widetilde \sigma_{k+1}$  such that  $\overline V_{22}$ is nonzero. Then for matrices $M, N$  given by Corollary 4.7, we have the relation
$$
K=MN=T_1G(x_n)-T_2,
$$
 where    $G(x_n)=[x_n^T\quad -1]\otimes I_{p+q}$, and
$$
\begin{array}l
  T_1={2\rho^{-2}{\cal K}xu^T}-[C_A^\dag \quad {\cal K}A^T],\qquad
T_2={\cal K}[I_n\quad 0_{n\times 1}]\otimes u^T,
\end{array}
$$
for  $ C_A^\dag=(I_n-{\cal K}A^TA)C^\dag$,  $\rho=\sqrt{1+\|x\|_2^2}$ and  $u^T=\big[-r^T(\widetilde A\widetilde C^\dag)\quad r^T\big]$ with  $r=Ax-B$. }

{Proof}. Note that when $k=n-p, d=1$,  in (\ref{4.7}) and (\ref{4.8}), $\Pi_{(d,n)}=\Pi_{(n,d)}=I_n$,   $\widetilde\Sigma_2^T\widetilde\Sigma_2=\widetilde\sigma_{n-p+1}^2$ and
$$
D^{-1}:=(  S_{1}^2\otimes I_{d}-\left[\begin{array}{cc} 0_p&0\\ 0& I_{n-p}\end{array}\right]\otimes ( \widetilde  \Sigma_{2}^T \widetilde \Sigma_{2}))^{-1}={\rm diag}(S_C^{-2}, (\widetilde\Sigma_1^2-\widetilde\sigma_{n-p+1}^2I_{n-p})^{-1}).
$$
 Following (\ref{4.00}), the matrices $M, N$ in (\ref{4.7})--(\ref{4.8})
 take the form
$$
\begin{array}l
M=\overline V_{22}^{-T}\widehat  V_{11}^{-T}D^{-1}=\overline V_{22}^{-T}[I_n\quad x_n]\widehat V_1D^{-1}=\overline V_{22}^{-T}[I_n\quad x_n]\widehat V_1{\rm diag}(S_C^{-2}, (\widetilde\Sigma_1^2-\widetilde\sigma_{n-p+1}^2I_{n-p})^{-1}),\\
N=N_1+N_2=[0_{(n+1)\times p}\quad \overline V_1]^T\otimes ({\bf Q}\widetilde U_2\widetilde \Sigma_2)^T+ \overline V_2^T\otimes
\Big([{\bf P}U_CS_C\quad {\bf Q}\widetilde U_1\widetilde \Sigma_1 ]^T \Big),
\end{array}
$$
where $\widehat V_1=[\widehat V_{11}^T\quad \widehat V_{21}^T]^T$ and $\widehat V_{11}$,  $\widehat V_{21}$ are defined in Theorem 4.2.

In the following, we first derive an equivalent formula for $(\widetilde\Sigma_1^2-\widetilde\sigma_{n-p+1}^2I_{n-p})^{-1}$. Let $\widetilde Q_2$ be given by (\ref{2.001}),
based on the SVD in (\ref{2.4}):  $\widetilde A\widetilde Q_2=\widetilde U_1\widetilde\Sigma_1\widetilde V_1^T+\widetilde U_2\widetilde\Sigma_2\widetilde V_2^T,$   partition
$$
\begin{array}{rl}
\widetilde V=[ \widetilde V_1\quad \widetilde V_2]&=\begin{array} c n-p\\ 1\end{array}\left[\begin{array}{cc} \widetilde V_{11}& \widetilde V_{12}\\
\widetilde V_{21}&\widetilde V_{22}\end{array}\right].\\
&\qquad\qquad  n-p\qquad 1
\end{array}
$$
Note that $AQ_2$ is the first $n-p$ columns of $\widetilde A\widetilde Q_2$, therefore $AQ_2=\widetilde U_1\widetilde\Sigma_1\widetilde V_{11}^T+\widetilde U_2\widetilde\Sigma_2\widetilde V_{12}^T$, from which
\begin{equation}
\begin{array}{rl}
Y_0&=(AQ_2)^T(AQ_2)-\widetilde \sigma_{n-p+1}^2I_{n-p}\\
&=[\widetilde V_{11}\quad \widetilde V_{12}]\Big({\rm diag}(\widetilde \Sigma_1^T\widetilde \Sigma_1, \widetilde \Sigma_2^T\widetilde \Sigma_2)-\widetilde \sigma_{n-p+1}^2I_{n-p+1}\Big)[\widetilde V_{11}\quad \widetilde V_{12}]^T\\
&=\widetilde V_{11}(\widetilde\Sigma_1^2-\widetilde\sigma_{n-p+1}^2I_{n-p})\widetilde V_{11}^T,
\end{array}\label{4.11}
\end{equation}
in which $\widetilde V_{1j}$  satisfies
\begin{equation}
\overline V_{1j}=Q_2\widetilde V_{1j}+\beta^{-1}x_C\widetilde V_{2j}=Q_2\widetilde V_{1j}-x_C\overline V_{2j},\quad j=1,2,\label{4.12}
\end{equation}
according to the  relation  $\overline V=\widetilde Q_2\widetilde V$.
Therefore  $\widetilde V_{22}=-\beta \overline V_{22} $   is  nonzero. By Lemma 4.1(a), $\widetilde V_{11}$ is nonsingular. From (\ref{4.11}), we obtain
$$
(\Sigma_1^2-\widetilde\sigma_{n-p+1}^2I_{n-p})^{-1}=\widetilde V_{11}^TY_0^{-1}\widetilde V_{11}.
$$

It should be noted that for any column vector $z$ and matrices $M_i$,
$$
M_1(M_2\otimes z^T)=(M_1M_2)\otimes z^T,\qquad z^T\otimes M_3=M_3(z^T\otimes I),
$$
therefore with $\widehat V_1=[V_C\quad \overline V_1]$,
$$
\begin{array}{rl}
MN_1&=\left([I_n\quad x_n]\widehat V_1\left[\begin{array}c
0_{p\times (n+1)}\\
\widetilde V_{11}^TY_0^{-1}\widetilde V_{11}\overline V_1^T\end{array}\right]\right)\otimes  ({\bf Q}\widetilde U_2\widetilde \Sigma_2\overline V_{22}^{-1})^T\\[5pt]
&=\big([I_n\quad x_n]\overline V_1\widetilde V_{11}^TY_0^{-1}\widetilde V_{11}\overline V_1^T\big)\otimes ({\bf Q}\widetilde U_2\widetilde \Sigma_2\overline V_{22}^{-1})^T,\\[5pt]
MN_2&=\widehat V_{11}^{-T}D^{-1} \Big([-x_n^T\quad 1]\otimes [{\bf P}U_CS_C\quad {\bf Q}\widetilde U_1\widetilde \Sigma_1 ]^T\Big)\\
&=\widehat V_{11}^{-T}D^{-1}[{\bf P}U_CS_C \quad {\bf Q}\widetilde U_1\widetilde \Sigma_1 ]^T \big([-x_n^T\quad 1]\otimes I_{p+q}\big)\\
&=[I_n\quad x_n]\widehat V_1[{\bf P}U_CS_C^{-1}\quad {\bf Q}\widetilde U_1\widetilde \Sigma_1\widetilde V_{11}^TY_0^{-T}\widetilde V_{11} ]^T ([-x_n^T\quad 1]\otimes I_{p+q})\\[5pt]
&=[I_n\quad x_n]\big(  V_CS_C^{-1}U_C^T{\bf P}^T+\overline V_1\widetilde V_{11}^TY_0^{-1}(\widetilde U_1\widetilde \Sigma_1 \widetilde V_{11})^T{\bf Q}^T)([-x_n^T\quad 1]\otimes I_{p+q})\\[5pt]
&=[I_n\quad x_n](\widetilde C^{\dag }{\bf P}^T+\overline V_1\widetilde V_{11}^TY_0^{-1}(\widetilde U_1\widetilde \Sigma_1 \widetilde V_{11})^T{\bf Q}^T)([-x_n^T\quad 1]\otimes I_{p+q}),
\end{array}
$$
in which
\begin{equation}
\begin{array}{rl}
\widetilde U_2\widetilde \Sigma_2\overline V_{22}^{-1}&=\widetilde A\widetilde Q_2\widetilde V_2\overline V_{22}^{-1}=\widetilde A\overline V_2\overline V_{22}^{-1}
=[A\quad b]\Big[{\displaystyle -x_n\atop\displaystyle 1}\Big]=:-r.\\
\widetilde U_1\widetilde \Sigma_1\widetilde V_{11}^T&=\widetilde U_1\widetilde \Sigma_1\widetilde V_{1}^T\Big[{\displaystyle  I_{n-p}\atop\displaystyle  0}\Big]
=[(\widetilde A\widetilde Q_2)-\widetilde U_2\widetilde\Sigma_2\widetilde V_2^T]\Big[{\displaystyle  I_{n-p}\atop\displaystyle  0}\Big]\\[5pt]
&=AQ_2-\widetilde U_2\widetilde\Sigma_2\overline V_{22}^{-1} \overline V_{22}\widetilde V_{12}^T=AQ_2+r\overline V_{22}\widetilde V_{12}^T,
\end{array}\label{4.13}
\end{equation}
and
\begin{equation}
\begin{array}{rl}
\widetilde V_{11}\overline V_1^T&=[\widetilde V_{11}\widetilde V_{11}^T\quad \widetilde V_{11}\widetilde V_{21}^T]\widetilde Q_2^T=[I_{n-p}-\widetilde V_{12}\widetilde V_{12}^T\quad
-\widetilde V_{12}\widetilde V_{22}^T]\widetilde Q_2^T\\
&=\Big([I_{n-p}\quad 0]-\widetilde V_{12}\widetilde V_2^T\Big)\widetilde Q_2^T=[Q_2^T\quad 0]-\widetilde V_{12}\overline V_2^T\\
&=[Q_2^T\quad 0]-\widetilde V_{12}\overline V_{22}^T[\overline V_{22}^{-T}\overline V_{12}^T\quad  1]=[Q_2^T\quad 0]+\widetilde V_{12}\overline V_{22}^T[x_n^T\quad -1].
\end{array}\label{4.14}
\end{equation}
Therefore
\begin{equation}
[I_n\quad x_n]\overline V_1\widetilde V_{11}^T=Q_2.\label{4.15}
\end{equation}
Moreover,  note that the Grevill's method \cite[Chapter 7, Section 5]{big} gives
$$
\widetilde C^\dag=\left[\begin{array}c \left(I_n- {\omega^{-1}{x_{\rm C}x_{\rm C}^T}}\right)C^\dag \\
{\omega^{-1} x_{\rm C}^TC^\dag }\end{array}\right],\qquad \omega= 1+\|x_{\rm C}\|_2^2.
$$
Combining this with
  the expression for $x_n$ in (\ref{2.002}) and  the residual $r_C=Ax_C-b$,   we have
\begin{equation}
\begin{array}l
[I_n\quad x_n]\widetilde C^\dag=(I_n- {\omega^{-1}{x_{\rm C}x_{\rm C}^T}}+\omega^{-1}x_nx_{C}^T)C^\dag=(I_n-\omega^{-1}{\cal K}A^Tr_Cx_C^T)C^\dag,\\
\widetilde A\widetilde C^\dag=(A-\omega^{-1}Ax_Cx_C^T+\omega^{-1}bx_C)C^\dag=(A-\omega^{-1}r_Cx_C^T)C^\dag.
\end{array}\label{4.16}
\end{equation}

 Combining (\ref{4.13})-(\ref{4.16}) with (\ref{4.12}), with $u^T=r^T{\bf Q}^T=\big[-r^T(\widetilde A\widetilde C^\dag)\quad r^T\big]$, ${\cal K}=Q_2Y_0^{-1}Q_2^T$,  we have
$$
\begin{array}{rl}
MN_1&=-\big[Q_2Y_0^{-1}Q_2^T\Big([I_{n-p}\quad 0]- Q_2\widetilde V_{12}\overline V_{22}^T[x_n^T\quad -1]\Big)\big]\otimes u^T\\[5pt]
&=-\big[(Q_2Y_0^{-1}Q_2^T)\Big([I_{n-p}\quad 0]+ (\overline V_{12}+x_C\overline V_{22})\overline V_{22}^T[x_n^T\quad -1]\Big)\big]\otimes u^T\\[5pt]
&=-\big[{\cal K}\Big([I_{n-p}\quad 0]+ (-x_n+x_C)\overline V_{22}\overline V_{22}^T[x_n^T\quad -1]\Big)\big]\otimes u^T\\[5pt]
&=-\Big([{\cal K}\quad 0]-\rho^{-2}{\cal K}x_n[x_n^T\quad -1]\Big)\otimes u^T\\[5pt]
&=-{\cal K}([I_n\quad 0]\otimes u^T)+(\rho^{-2}{\cal K}x_nu^T)([x_n^T\quad -1]\big)\otimes I_{p+q}),\\[5pt]
MN_2&=\Big([(I_n-\omega^{-1}{\cal K}A^Tr_Cx_C^T)C^\dag\quad 0]+Q_2Y_0^{-1}Q_2^TA^T[-\widetilde A\widetilde C^\dag\quad I_q]\\[5pt]
&\quad +(Q_2Y_0^{-1}Q_2^T) Q_2\widetilde V_{12}\overline V_{22}u^T\Big)([-x_n^T\quad 1]\otimes I_{p+q})\\[5pt]
&=\Big([(I_n-\omega^{-1}{\cal K}A^Tr_Cx_C^T)C^\dag\quad 0]+[-{\cal K}A^T(A-\omega^{-1}r_Cx_C^T)C^\dag\quad {\cal K}A^T]\\[5pt]
&\quad +{\cal K}(\overline V_{12}+x_C\overline V_{22})\overline V_{22}^Tu^T\Big)([-x_n^T\quad 1]\otimes I_{p+q})\\[5pt]
&=\Big([(I_n-{\cal K}A^TA)C^\dag\quad {\cal K}A^T]+{\cal K}(-x_n+x_C)\overline V_{22}\overline V_{22}^Tu^T\Big)([-x_n^T\quad 1]\otimes I_{p+q})\\[5pt]
&=-\Big([(I_n-{\cal K}A^TA)C^\dag\quad {\cal K}A^T]-\rho^{-2}{\cal K}x_nu^T\Big)\big([x_n^T\quad -1]\otimes I_{p+q}\big),
\end{array}
$$
where $\overline V_{22}\overline V_{22}^T=\|\overline V_{22}\|_2^2=\rho^{-2}$  based on the fact that
$[x_n^T\quad -1]=\rho\big[\overline V_{12}^T\quad \overline V_{22}^T\big]^T$ for $\rho^2= 1+\|x_n\|_2^2$, and ${\cal K}(x_C-x_n)=-{\cal K}x_n$ since $Q_2^Tx_C=0$. Therefore
$$
M(N_1+N_2)=\Big(2\rho^{-2}{\cal K}x_nu^T-[(I_n-{\cal K}A^TA)C^\dag\quad {\cal K}A^T]\Big)\big([x_n^T\quad -1]\otimes I_{p+q}\big)-{\cal K}([I_n\quad 0]\otimes u^T\big),
$$
which is exactly $K$. The assertion in the theorem then follows.\qed\\

{\bf Remark 4.1 }  When $p>0, k=n-p$ and $d=1$,  for the single dimensional TLSE problem,  Liu and Jia \cite{liu2} derived the first order perturbation estimate as
$$
\Delta x=K_{L, H}{\rm vec}([\Delta L\quad \Delta H])+{\cal O}(\|[\Delta L\quad \Delta H]\|_F^2),
$$
where $K_{L,H}=MN=K$. With this,   three types of  condition number formulae of the single-dimensional TLSE problem were derived.
 The result in Corollary 4.8 shows that the newly derived perturbation analysis and condition numbers for the multidimensional case unify those for the single dimensional TLSE problem.

 .

{\bf Remark 4.2} It is observed that the formulae for three types of condition numbers involve the Kronecker product which might lead to large storage and computation cost.
For mixed and componentwise condition numbers, we can use their upper bounds as alternations, while for the normwise condition numbers, as did in \cite{zmw}, we can compute
$$
\kappa^{\rm abs}=\|{\breve H}\|_2=\|{\breve H}^T{\breve H}\|_2^{1/2},\quad for\quad {\breve H}=(H_1+H_2)G\overline Z=(H_1+H_2)D^{-1}{\breve Z},
$$
by applying the power method to the matrix ${\breve H}^T{\breve H}$, in which $D$ is defined in Theorem 4.6, and $\breve Z=\big[ I_t\otimes \widetilde \Sigma_{2}^T\qquad S_{1}\otimes I_{n+d-t}\big]\overline Z$. In the power scheme,
  the matrix-vector multiplications associated  with ${\breve H}$ and ${\breve H}^T$ can be transformed into Kronecker product-free operations, say for ${\breve H}f$, where $f=[f_1^T\quad f_2^T]^T$ with $f_i={\rm vec}(F_i)$ with $F_1\in {\mathbb R}^{t\times (p+q)}, F_2\in {\mathbb  R}^{(n+d-t)\times t}$,
  $$
  \begin{array}l
  g:={\breve H}f=(H_1+H_2)D^{-1}{\rm vec}\Big(({\bf Q}\widetilde U_2\widetilde\Sigma_2)^TF_1\left[\begin{array}{cc} 0_p&0\\ 0& I_k\end{array}\right] +F_2\left[\begin{array}{cc}
  S_C&-U_C^T(\widetilde A\widetilde C^\dag)^T\widetilde U_1\widetilde \Sigma_1\\ 0&\widetilde\Sigma_1\end{array}\right]  \Big)\\
  =(H_1+H_2){\rm vec}(T)
  ={\rm vec}\Big((\overline V_{12}+X_t\overline V_{22}))T\widehat V_{21}^T(\overline V_{22}\overline V_{22}^T)^{-1}  +(\overline V_{11}+X_t\overline V_{21})T^T\overline V_{22}^\dag\Big),
  \end{array}
  $$
  where  $t_i=Te_i$ satisfies
  $$
  t_i=(s_i^2I_{n+d-t}-\tau_i\widetilde\Sigma_2^T\widetilde\Sigma_2)^{-1}\Big(({\bf Q}\widetilde U_2\widetilde\Sigma_2)^TF_1\left[\begin{array}{cc} 0_p&0\\ 0& I_k\end{array}\right] +F_2\left[\begin{array}{cc}
  S_C&-U_C^T(\widetilde A\widetilde C^\dag)^T\widetilde U_1\widetilde \Sigma_1\\ 0&\widetilde\Sigma_1\end{array}\right]  \Big)e_i,
  $$
  in which $s_i, \tau_i$ are the same as those in Theorem 4.6.
The Kronecker product-free expression associated with ${\breve H}^Tg$ can be derived in a similar manner. Here we omit these.

\section{Numerical experiments}
\numberwithin{table}{section}

In this section, we present numerical examples to verify   our results. The following numerical tests are performed via MATLAB
  with machine precision $u = 2.22e-16$ in a laptop with Intel Core (TM) i5-5200U CPU by using double precision.\\

{\it Example 1} In this example, we generate  random multidimensional TLSE problems to verify the rationality of the first order perturbation estimate in Theorem 4.2.
The entries in $[C\quad D]$  and $[A\quad B]$ are generated as random variables uniformly distributed in the interval (0,1), via Matlab command `{\sf rand}($\cdot$)'.
Set $p=10, q=40, n=40, d=5$, and let the perturbations to the data be given by
$$
[\Delta C\quad \Delta D]=\epsilon* {\sf rand}(p, n+d),\qquad [\Delta A\quad \Delta B]=\epsilon* {\sf rand}(q, n+d).
$$
Choose $t=10, 20, 30, 40$  and compute the solutions to the original and perturbed problems via the QR-SVD method. In Table \ref{tab5.1} we compute the absolute error
$$
\eta_{\Delta X_t}=\left\| {\rm vec}(\Delta X_t)-K{\rm vec}([\Delta L\quad \Delta H])\right\|_\infty,
$$
with respect to different $\epsilon$.

 \renewcommand\tabcolsep{22.0pt}
\begin{table}
\begin{center}
\caption{The absolute error of the first order perturbation estimate of ${\rm vec}(\Delta X_{t})$}
\begin{tabular}{ccccc}
\hline
$t$&$10$&$20$&$30$&$40$\\\hline
$\epsilon=10^{-2}$&1.9e-4&6.3e-4&5.2e-4&3.0e-4\\
$\epsilon=10^{-4}$&5.6e-8&2.3e-8&3.7e-8&2.1e-8\\
$\epsilon=10^{-6}$&2.7e-12&3.7e-12&1.8e-12&1.5e-12\\\hline
\end{tabular}\label{tab5.1}
\end{center}
\end{table}

The tabulated results show that $\eta_{\Delta X_t}={\cal O}(\epsilon^2)$, illustrating the rationality of the first order perturbation estimates in Theorem 4.2.

{\it Example 2.} In this example, we do some numerical experiments for TLSE from piecewise-polynomial data fitting problem that is modified from \cite[Chapter 16]{bv} and \cite{liu1}.

 Given $N$ points $(t_i, y_i)$ on the plane, we are seeking to find a piecewise-polynomial function $f(t)$ fitting the above set of the points, where
 $$
 f(t)=\left\{\begin{array}{ll}
 f_1(t),&t\le a,\\
 f_2(t),&t>a,
 \end{array}\right.
 $$
 with $a$ given, and $f_1(t)$ and $f_2(t)$ polynomials of degree three or less,
 $$
 f_1(t)=x_1+x_2t+x_3t^2+x_4t^3,\qquad  f_2(t)=x_5+x_6t+x_7t^2+x_8t^3,.
 $$
 The conditions that $f_1(a)=f_2(a)$ and $f_1'(a)=f_2'(a)$ are imposed, so that $f(t)$ is continuous and has  a continuous first derivative at $t=a$. Suppose   the $N$ data are numbered so that $t_1,\ldots,t_M\le a$ and $t_{M+1},\ldots, t_N>a$. The conditions $f_1(a)-f_2(a)=0$ and $f_1'(a)-f_2'(a)=0$ leads to the equality constraint $Cx=d$ for $x=[x_1,x_2,\ldots,x_8]^T$ and
$$
C=\left[\begin{array}{cccccccc}
1&a&a^2&a^3&-1&-a&-a^2&-a^3\\
0&1&2a&3a^2&0&-1&-2a&-3a^2\end{array}\right],\quad d=\left[\begin{array}c 0\\0\end{array}\right].
$$
The vector $x$ that minimizes the  sum of squares of the prediction errors
 $$
 \sum\limits_{i=1}^M(f_1(t_i)-y_i)^2+\sum\limits_{i=M+1}^N(f_2(t_i)-y_i)^2,
 $$
 gives $\min_x\|Ax- b\|_2$, where
 $$
 A=\left[\begin{array}{cccccccc}
 1&t_1&t_1^2&t_1^3&0&0&0&0\\
 1&t_2&t_2^2&t_2^3&0&0&0&0\\
 \vdots&\vdots&\vdots&\vdots& \vdots&\vdots&\vdots&\vdots\\
 1&t_M&t_M^2&t_M^3&0&0&0&0\\
 0&0&0&0&1&t_{M+1}&t_{M+1}^2&t_{M+1}^3\\
 0&0&0&0&1&t_{M+2}&t_{M+2}^2&t_{M+2}^3\\
 \vdots&\vdots&\vdots&\vdots& \vdots&\vdots&\vdots&\vdots\\
  0&0&0&0&1&t_{N}&t_{N}^2&t_{N}^3
  \end{array}\right],\qquad b=\left[\begin{array}c y_1\\y_2\\\vdots\\y_M\\y_{M+1}\\\vdots\\ y_N\end{array}\right],
 $$
 and the matrix $A$ is of 50\% sparsity. If more than one observation vector is allowed, the data fitting problem becomes the multidimensional TLSE problem (\ref{1.2}).

Take $M=200, N=200$  and  sample $t_i\in [0, 1]$ randomly. For a randomly generated piecewise-polynomial function $f (t)$ with a predetermined $a$, we
compute the corresponding function value $y_i=f(t_i)$. We add  random componentwise perturbations on the data
as
\begin{equation}
\Delta L=10^{-12}\cdot E_{p+q,8}\odot L, \qquad \Delta H=10^{-12}\cdot E_{p+q,d}\odot H,\quad q=M+N,\label{5.1}
\end{equation}
 where  $E_{s,t}$ is the random $s\times t$ matrix whose entries are uniformly distributed on the interval (0,1), $\odot$ denotes the entrywise multiplication.

 For simplicity let $\kappa_{\rm n}, m, c$ denote the relative normwise condition number, mixed and componentwise condition numbers given in
Theorem 4.3 and Theorem 4.6, respectively.  Set
$$
x={\rm vec}(X_t),\quad \epsilon_{\rm n}={\|[\Delta L\quad \Delta H]\|_F\over \|[L\quad H]\|_F},\qquad
\epsilon_{\rm c}=\min \{\epsilon: |\Delta L|\le \epsilon |L|,  |\Delta H|\le \epsilon |H|\},
$$
where $t$ is a random integer between $p$ and $n$ such that $\overline V_{22}$ is of full row rank,  and the quantity  $\epsilon_{\rm n}$     is used to evaluate the upper bound of the  forward error ${\|\Delta x\|_2\over \|x\|_2}$  via $\epsilon_{\rm n}\kappa_{\rm n}$,
while $\epsilon_{\rm c}$ is to derive upper bounds for ${\|\Delta x\|_\infty\over \|x\|_\infty}, \|{\Delta x\over x}\|_\infty$  via  mixed and componentwise condition numbers.
Moreover we  let  $\rho=\rho_{AC}^{(2)}\eta_k^\sigma$ be the factor of upper bounds  of $\kappa^{\rm abs}(X_t, L, H)$.

We list numerical results with respect to different $a$, and for each $a$, we generate two different problems and compare the estimated upper bound with actual relative  forward errors. It's observed that for fixed $a$, the problems with a larger $\|X_t\|_2^2$ and moderate $\rho$  produce larger condition number estimates, which illustrates that  the norm $\|X_t\|_2$ is a factor to affect the condition number of TLSE problem.
However whether $\|X_t\|_2$ is big or small, the  estimated upper bounds of the forward error via $\epsilon_{\rm n}\kappa_{\rm n}$, $\epsilon_{\rm c}m$, $\epsilon_{\rm c}c$ are about one or two orders of magnitude larger than the corresponding forward error of the solution.
Among three upper bounds $\kappa_{\rm n}^{u}, m^u, c^u$ of condition numbers,  the normwise condition number-based upper bound $\kappa_{\rm n}^{u}$ is acceptable and is
about one or two orders of magnitude larger than $\kappa_{\rm n}$. The upper bounds $m^u, c^u$ are more sharp, which are  at most one order of magnitude larger than the corresponding exact condition numbers, therefore they are good estimates of their corresponding condition numbers and forward error of the solutions.

\renewcommand\tabcolsep{3.0pt}
\begin{table}[!htb]
\begin{center}
\caption{Comparisons of forward errors and upper bounds for  the perturbed TLSE problem}
\label{Table1}
\begin{tabular}{*{28}{l}}
\hline
$a$&$\|X_t\|_2^2$&$\rho$&$\frac{\|\Delta x\|_2}{\|x\|_2}$&$\epsilon_{\rm n}\kappa_{\rm n}$&$\epsilon_{\rm n}\kappa_{\rm n}^{u}$&   ${\|\Delta x\|_\infty\over \|x\|_\infty}$ & $\epsilon_{\rm c}m$&
$\epsilon_{\rm c}m^u$&$\|{\Delta x \over x}\|_\infty$ & $\epsilon_{\rm c}c$&$\epsilon_{\rm c}c^u$\\
 \hline

0.1&4.2& 12.0& 2.2e-13& 2.1e-11& 7.7e-10& 2.6e-13& 3.8e-12& 1.2e-11& 7.3e-13& 6.9e-12& 3.2e-11\\
& 1.6e5& 76.0& 1.1e-11& 3.0e-9& 7.7e-7& 9.8e-12& 9.1e-10& 1.7e-9& 1.5e-10& 2.1e-8& 2.8e-8\\\hline

0.3& 4.2& 12.0& 1.2e-13& 2.1e-11& 7.5e-10& 1.7e-13& 4.0e-12& 1.4e-11& 5.2e-13& 7.2e-12& 3.1e-11\\
 &2.4e5& 42.0& 7.4e-12& 3.5e-9& 5.2e-7& 6.4e-12& 1.2e-9& 3.0e-9& 4.3e-11& 6.4e-9& 1.3e-8\\\hline

0.5& 5.6& 12.0& 1.8e-13& 4.7e-11& 8.2e-10& 2.3e-13& 1.2e-11& 3.3e-11& 2.7e-11& 2.0e-9& 7.2e-9\\
 &5.3e4& 68.0& 5.8e-12& 1.0e-9& 3.8e-7& 6.8e-12& 7.8e-10& 1.6e-9& 1.2e-9& 1.0e-7& 1.6e-7\\\hline

0.7& 3.0& 11.0& 1.2e-13& 2.2e-11& 7.1e-10& 1.4e-13& 4.8e-12& 1.6e-11& 2.7e-13& 8.4e-12& 3.0e-11\\
&1.3e7& 75.0& 1.4e-10& 1.5e-8& 7.0e-6& 1.4e-10& 6.5e-9& 1.1e-8& 2.3e-8& 1.1e-6& 1.9e-6\\\hline

0.9& 2.3& 11.0& 4.0e-14& 2.2e-11& 6.9e-10& 5.8e-14& 5.7e-12& 1.9e-11& 8.0e-14& 7.9e-12& 2.8e-11\\
 & 5.2e8& 42.0& 1.0e-9&1.0e-7& 2.7e-5& 1.3e-9& 7.6e-8& 2.1e-7& 1.7e-9& 5.4e-7& 1.0e-6\\\hline

\end{tabular}\label{tab5.2}
\end{center}
\end{table}

{\it Example 3.} This example is modified from \cite{bg}. Let $p=d=5, n=10,q=20, k=3, t=p+k=8$, and
 $\widetilde Q$ be an  arbitrary $(n+d)\times (n+d)$ orthogonal matrix  and $\widetilde Q_1$  is the submatrix of $\widetilde Q$ by taking its first $p$ columns.
Let   $U_0$ be an arbitrary $p\times p$ orthogonal matrix,  $y, z$ be unit column vectors of length $q, n+d$, respectively,  set
$$
\begin{array}l
\widetilde C=[C~~ D]=U_0{\rm diag}([1, 0.5, 0.1, 0.1, \kappa_C^{-1}])\widetilde Q_1^T,\quad \widetilde A=[A~B]=\hat A\widetilde Q^T,\quad \mbox{with}\\
\hat A=(I_q-2yy^T)[\hat \Sigma\quad O](I_{n+d}-2zz^T),\\
\hat \Sigma={\rm diag}(10, 8, 1, 1, 1, 1, 1, 1-\delta/2,  1-\delta, 1-2\delta,  1/6, 1/7,\cdots, 1/10]),
\end{array}
$$
where $\kappa_C$  is used to control the condition number of $[C~ D]$. Note that $\widetilde A\widetilde Q_2$ is the last $n+d-p$ columns of $\widetilde A\widetilde Q$, and
 by the interlacing theorem of the singular values, the relation $1=\sigma_j(\widetilde A\widetilde Q)\ge \sigma_j(\widetilde A\widetilde Q_2)\ge \sigma_{p+j}(\widetilde A\widetilde Q),$  for $j=k, k+1$ and therefore $0<\delta<1/12$ can be  used to control the gap of the singular values $\widetilde \sigma_k, \widetilde \sigma_{k+1}$ of $\widetilde A\widetilde Q_2$.

Consider the same perturbation as in (\ref{5.1}),  for different $\kappa_C$ and $\delta$,  we compute the forward errors and upper bounds via three types condition numbers in Table \ref{tab5.3}.
It's observed that the  estimated upper bounds of the forward error via $\epsilon_{\rm n}\kappa_{\rm n}$, $\epsilon_{\rm c}m$, $\epsilon_{\rm c}c$ are about one or two orders of magnitude larger than the corresponding forward error of the solution, even the quantity $\rho$ is very large.
For the compact upper bounds  $ m^u, c^u$ of condition numbers,  $m^u, c^u$ are very sharp
in most cases, while   $\kappa_{\rm n}^{u}$  is not robust against the ill-conditioning of $\widetilde C$ and sometimes they are  three  orders of magnitude larger than $\kappa_{\rm n}$ and five or six  orders of magnitude larger than $\frac{\|\Delta x\|_2}{\|x\|_2}$.

\renewcommand\tabcolsep{3.0pt}
\begin{table}[!htb]
\begin{center}
\caption{Comparisons of forward error and upper bounds for  the perturbed TLSE problem}
\label{Table1}
\begin{tabular}{*{28}{l}}
\hline
$\sigma$&$\|X_t\|_2^2$&$\rho$&$\frac{\|\Delta x\|_2}{\|x\|_2}$&$\epsilon_{\rm n}\kappa_{\rm n}$&$\epsilon_{\rm n}\kappa_{\rm n}^{u}$&   ${\|\Delta x\|_\infty\over \|x\|_\infty}$ & $\epsilon_{\rm c}m$&
$\epsilon_{\rm c}m^u$&$\|{\Delta x \over x}\|_\infty$ & $\epsilon_{\rm c}c$&$\epsilon_{\rm c}c^u$\\\hline

\multicolumn{12}{l}{$\kappa_C=10^1$}\\\hline
0.1&2.1& 2.4e3& 2.1e-12& 5.4e-10& 3.9e-8& 1.7e-12& 3.4e-11& 7.8e-11& 3.2e-11& 8.6e-10& 2.1e-9\\
0.01& 0.71& 6.7e3& 4.2e-12& 2.0e-10& 6.5e-8& 4.9e-12& 1.4e-10& 1.5e-10& 1.1e-8& 1.9e-7& 3.0e-7\\
 0.001&0.9& 1.6e5& 1.1e-10& 1.2e-8& 2.4e-6& 1.4e-10& 3.1e-9& 4.2e-9& 3.6e-10& 1.1e-8& 1.3e-8\\\hline

\multicolumn{12}{l}{$\kappa_C=10^3$}\\\hline
0.1&0.74& 1.2e5& 1.1e-10& 1.5e-8& 1.4e-6& 1.2e-10& 1.1e-9& 2.3e-9& 6.5e-10& 1.4e-8& 4.5e-8\\
0.01& 2.1& 9.8e5& 7.7e-11& 8.2e-9& 8.7e-6& 6.7e-11& 6.3e-10& 1.4e-9& 4.7e-9& 3.6e-8& 9.5e-8\\
0.001&0.51& 5.5e6& 1.2e-10& 2.7e-8& 7.7e-5& 1.5e-10& 3.3e-9& 5.3e-9& 3.6e-9& 9.9e-8& 2.4e-7\\\hline

\multicolumn{12}{l}{$\kappa_C=10^6$}\\\hline
0.1&2.4& 8.9e7& 2.7e-8& 9.7e-6& 1.4e-3& 3.4e-8& 6.2e-7& 1.2e-6& 6.6e-7& 2.7e-5& 1.3e-4\\
 0.01&4.8& 6.1e8& 1.2e-7& 3.5e-5& 1.2e-2& 1.1e-7& 1.2e-6& 2.9e-6& 2.2e-6& 2.8e-5& 6.7e-5\\
 0.001&2.0& 4.6e9& 2.0e-8& 1.1e-5& 6.5e-2& 2.1e-8& 2.5e-7& 4.6e-7& 4.6e-7& 5.5e-6& 1.8e-5\\\hline
\end{tabular}\label{tab5.3}
\end{center}
\end{table}

\section{Conclusion}
In this paper, we investigate the  solution of multidimensional TLSE problem, and prove it is equivalent to the multidimensional weighted TLS solution in the limit sense, with the aid of perturbation theory of invariant subspace. Based on this close relation, the closed formula for the first order perturbation estimate of the minimum Frobenius norm TLSE solution
 $X_t=-\overline V_{12}\overline V_{22}^\dag$  is derived, from which the expressions for normwise, mixed and componentwise condition numbers of problem  TLSE are also presented.
  Since there expressions involve matrix Kronecker product operations which may make the computation more expensive, we provide compact upper bounds to
  enhance the computation efficiency. All  expressions and upper bounds of these condition numbers  generalize those for the single-dimensional TLSE problem \cite{liu2} and multidimensional TLS problem \cite{mzw}.

Some numerical examples are also given in this paper to demonstrate the effectiveness in estimating the forward errors.
Tightness of  upper bounds for mixed and componentwise condition numbers are   shown in numerical examples, even for ill-conditioned problems, while it is not necessarily
true for the upper bounds of  the normwise  condition number.
Therefore in order to derive good estimates of forward errors via normwise condition number,    we recommend using power scheme to compute the true  value  to avoid  Kronecker product operations.


\begin{thebibliography}{}
%
\bibitem{bg} M. Baboulin,  S. Gratton, A contribution to the conditioning of the total least-squares problem, SIAM
J Matrix Anal. Appl., 32(3)  (2011), pp. 685-699.


\bibitem{big} A. Ben-Israel, T. N.E. Greville, Generalized inverses, theory and applications, 2nd ed., Spring-Verlag New York, (2003).


\bibitem{bv} S. Boyd, L. Vandenberghe, {Introduction to applied linear algebra-vectors, matrices, and least squares},  https://web.stanford.edu/~boyd/vmls/vmls.pdf, (2017)

\bibitem{ch} A. J. Cox, N. J. Higham,  {  Accuracy and stability of the null space method for solving the equality constrained least squares problem}.
  BIT, 39(1)(1999), pp. 34-50.



\bibitem{diao3} H. Diao,
 {  Condition numbers for a linear function of the solution of the
linear least squares problem with equality constraints}.
 Journal of Computational and Applied Mathematics,  344(2018),pp. 640-65.



\bibitem{ds} H. Diao, Y. Sun,  Mixed and componentwise condition numbers for a linear function
of the solution of the total least squares problem, Linear Algebra and its Applications
544:1(2018), pp. 1-29.



\bibitem{dwx}  H. Diao, Y. Wei, P. Xie,  Small sample statistical condition estimation for the total least squares
problem, Numer. Algorithms, {75}(2)  (2017), pp. 1-21.






\bibitem{ddl} E.M. Dowling, R.D. Degroat, D.A. Linebarger, {  Total least squares with linear constraints},  IEEE International Conference on Acoustics, 5(5), (1992), pp. 341-344.

\bibitem{ge} A.J. Geurts,  A contribution to the theory of condition, Numer Math., 39(1)  (1982), pp. 85-96.

\bibitem{gk}
I.  Gohberg,  I. Koltracht: Mixed, componentwise, and structured condition numbers. SIAM J. Matrix Anal. Appl.  {  14} (1993), pp. 688-704


\bibitem{gti} S. Gratton, D. Titley-Peloquin,  J. T. Ilunga, Sensitivity and conditioning of
the truncated total least squares solution, SIAM J. Matrix Anal. Appl., 34 (2013), pp. 1257-1276.


\bibitem{gv} G.H. Golub, C.F.  Van Loan,  {   An analysis of total least squares problem},  SIAM J Matrix Anal Appl., 17(6) (1980), pp. 883-893.

\bibitem{gv2} G.H. Golub, C.F. Van Loan, { Matrix Computations(4ed.)},
 Johns Hopkins University Press, Baltimore (2013)

 \bibitem{gr}
A. Graham,   Kronecker Products and Matrix Calculus with Application, Wiley, New York,
 MR0640865 (83g:15001) (1981)



\bibitem{hv} K. Hermus, W. Verhelst, P. Lemmerling, P. Wambacq,  S. Van Huffel, {  Perceptual audio modeling with
exponentially damped sinusoids},  Signal Processing, 85 (2005),  pp. 163-176.

 \bibitem{jl}
Z. Jia,  B. Li,  On the condition number of the total least squares problem. Numer. Math. {  125}(1) (2013), pp. 61-87.

\bibitem{ls}
 A.N. Langville,  W.J. Stewart,  The Kronecker product and stochastic automata
networks. J. Comput. Appl. Math. {  167}(2004), pp. 429-447.

\bibitem{ld} P. Lemmerling, B. De Moor, { Misfit versus latency},  Automatica, 37(2001), pp. 2057-2067.


\bibitem{lmv} P. Lemmerling, N. Mastronardi,  S. Van Huffel, { Efficient implementation of a structured total least squares
based speech compression method},  Linear Algebra Appl., 366(2003),  pp. 295-315.

\bibitem{lj}
B. Li, Z. Jia, Some results on condition numbers of the scaled total least squares problem. Linear
Algebra Appl. { 435}(3)  (2011), pp.  674-686.


\bibitem{liu1} Q. Liu, C. Chen, Q. Zhang, Perturbation analysis for total least squares problems with linear equality constraint, Applied Numerical Mathematics,
161(2021), pp. 69-81.

\bibitem{liu2} Q. Liu, Z. Jia,
{On   condition numbers of the total least squares problem with linear equality constraint},
 arxiv:2008.08233 [math.NA].


\bibitem{liu} Q. Liu, S. Jin, L. Yao, D. Shen,
The revisited total least squares problems with linear equality constraint, Applied numerical mathematics, 152(2020), pp. 275-284.





\bibitem{mdb} Q. Meng, H. Diao, Z. Bai,  Condition numbers for the truncated total least squares problem and their estimations, arXiv:2004.12082[math.NA]


\bibitem{mv} I. Markovsky, S. Van Huffel, { Overview of total least squares methods},  Signal Processing, 87(2007), pp. 2283-2302.



\bibitem{mzw} L. Meng,  B. Zheng and Y. Wei,
Condition numbers of the multi-dimensional total least squares problems
having more than one solution, Numerical Algorithms, 84 (2020) 887-908.



\bibitem{nbk} M. Ng, N. Bose, J. Koo,  { Constrained total least squares for color image reconstruction},
  Total Least Squares and Errors-in-Variables Modelling III: Analysis, Algorithms and Applications, Kluer Academic Publishers, (2002), pp. 365-374.

\bibitem{npp} M. Ng, R. Plemmons, F. Pimentel, {A new approach to constrained total least squares image restoration},
Linear Algebra Appls., 316(2000), pp. 237-258.




\bibitem{pe} K. Pearson, { On lines and planes of closest fit to systems of points in space},  Phil. Mag.,  2(1901), pp. 559-572.




\bibitem{pw}  A. R. De Pieero, M. Wei,  Some new properties of the equality constrained and weighted least squares problem, Linear Algebra and its Applications
{\bf 320}(1-3) (2000), pp. 145-165.




\bibitem{sc} B. Schaffrin, { A note on constrained total least squares estimation},
 Linear Algebra Appl. 417(2006), pp.245-258.

\bibitem{rice} J.R.  Rice, A theory of condition. SIAM J. Numer. Anal. 1966; {3}:287-310.

\bibitem{st} G.W. Stewart, On the asymptotic behavior of scaled singular value and QR decompositions, Mathematics of  Computation, 43(168)(1984), pp. 483-489.

\bibitem{sun} G. W. Stewart and J.-G. Sun, Matrix Perturbation Theory, Academic Press, Boston, 1990.

\bibitem{vv} S. Van Huffel, On the significance of nongeneric total least squares problems, SIAM J. Matrix
Anal. Appl., 13 (1992), pp. 20-35.


\bibitem{val} S. Van Huffel, P. Lemmerling, eds.  {Total Least Squares and Errors-in-Variables Modeling: Analysis, Algorithms and Applications},
  Kluwer, Dordrecht, Boston, London, (2002).

\bibitem{vh1} S. Van Huffel, J. Vandevalle,
{  The Total Least Squares Problems: Computational Aspects and Analysis},  Vol. 9 of Frontiers in Applied Mathematics, SIAM, Philadelphia, (1991).

\bibitem{vv2} S. Van Huffel and J. Vandewalle, Analysis and solution of the nongeneric total least squares
problem, SIAM J. Matrix Anal. Appl., 9 (1988), pp. 360-372.

\bibitem{wei} M. Wei, Algebraic relations between the total least squares and least squares problems with
more than one solution, Numer. Math., 62 (1992), pp. 123-148.

\bibitem{wei2} M. Wei, Perturbation theory for the rank-deficient equality constrained least squares problem. SIAM J. Numer. Anal. 29:5 (1992), pp. 1462-1481.

\bibitem{wp}   M. Wei,   A. R. De Pieero,  Upper perturbation bounds of weighted projections, weighted and constrained least squares
problems, SIAM J. Matrix Anal. Appl. {\bf 21}(3) (2000), pp. 931-951.


\bibitem{xxw} P. Xie, H. Xiang, Y. Wei, A contribution to perturbation analysis for total least
squares problems, Numerical Algorithms, 75(2) (2017), pp. 381-395.



\bibitem{xxw2} P. Xie, H. Xiang and Y. Wei, Randomized algorithms for total least squares problems,
Numer Linear Algebra Appl., 26 (2019) e2219.






\bibitem{zlw} L. Zhou, L. Lin, Y. Wei, S. Qiao,  Perturbation analysis and condition numbers of scaled total least
squares problems. Numer. Algorithms 51(3)(2009), pp. 381-399.



\bibitem{zmw} B. Zheng, L. Meng and Y. Wei, Condition numbers of the multidimensional total least squares problem.
SIAM J. Matrix Anal. Appl., 38 (2017), pp. 924-948.

\bibitem{zy} B.  Zheng ,  Z. Yang,  Perturbation analysis for mixed least squares-total least squares
problems. Numer Linear Algebra Appl. 2019;26:e2239. https://doi.org/10.1002/nla.2239


\end{thebibliography}


\end{document}